\documentclass{article}
\usepackage[utf8]{inputenc}
\usepackage[numbers]{natbib}
\usepackage{manfnt} 
\newcommand{\dnote}[1]{%
    \noindent % I guess this is intended...
    \begin{tabular}{@{}m{0.13\textwidth}@{}m{0.87\textwidth}@{}}%
        \huge\textdbend &#1%
    \end{tabular}%
    \par % ... and this too.
}
\usepackage{amsthm}

\newtheorem{theorem}{Theorem}
\newtheorem{lemma}{Lemma}
\newtheorem{proposition}{Proposition}

\newtheorem{definition}{Definition}
\usepackage[top=0.8in, bottom=1in, left=1.25in, right=1.25in]{geometry}
\usepackage{scrextend}
\usepackage[utf8]{inputenc} % allow utf-8 input
\usepackage[T1]{fontenc}    % use 8-bit T1 fonts
\usepackage[english]{babel}
\usepackage{lmodern}
\usepackage[bitstream-charter]{mathdesign}
\usepackage[table]{xcolor}
 \usepackage{hyperref}
 \definecolor{burgundy}{rgb}{0.5, 0.0, 0.13}
\definecolor{camel}{rgb}{0.76, 0.6, 0.42}
\definecolor{chamoisee}{rgb}{0.63, 0.47, 0.35}
\definecolor{grey1}{RGB}{128,128,128}
 \hypersetup{colorlinks=true,linkcolor=burgundy,citecolor=chamoisee,urlcolor=burgundy,linktoc=page}
\usepackage{url}            % simple URL typesetting
\usepackage{booktabs}       % professional-quality tables
\usepackage{nicefrac}       % compact symbols for 1/2, etc.
\usepackage{microtype}      % microtypography
\usepackage{dsfont}         % ds font
\renewcommand{\epsilon}{\varepsilon}
\renewcommand{\phi}{\varphi}

\title{Concentration inequality for U-statistics of order two for uniformly ergodic Markov chains}   
\author{%
  Quentin Duchemin\thanks{This work was supported by grants from Région Ile-de-France.} \\
 LAMA, Univ Gustave Eiffel, CNRS, Marne-la-Vallée, France.\\
  \texttt{quentin.duchemin@univ-eiffel.fr} \\
  $\And$\\
  Yohann De Castro\\
  Institut Camille Jordan, École Centrale de Lyon, Lyon, France \\
  \texttt{yohann.de-castro@ec-lyon.fr} \\
    $\And$\\
  Claire Lacour\\
  LAMA, Univ Gustave Eiffel, CNRS, Marne-la-Vallée, France. \\
  \texttt{claire.lacour@univ-eiffel.fr} 
}

\date{}

\usepackage{amsmath}
\usepackage{mathtools}

\DeclarePairedDelimiter\floor{\lfloor}{\rfloor}
\usepackage{diagbox}
\usepackage{subcaption}

\usepackage{algorithm,algorithmic,refcount}
\usepackage{wrapfig}
\usepackage{tikz}
\usepackage{makecell}

\usepackage{tikz}
\usetikzlibrary{%
  arrows,%
  calc,%
  shapes.geometric,%
  shapes.misc,%
  shapes.symbols,%
  shapes.arrows,%
  automata,%
  through,%
  positioning,%
  scopes,%
  decorations.shapes,%
  decorations.text,%
  decorations.pathmorphing,%
  shadows}
\usetikzlibrary{positioning}
\usepackage{multirow}
\usepackage{makecell}
\usepackage{mdframed}
\newtheorem{assumption}{Assumption}
\usepackage{savesym}
\savesymbol{checkmark}
\usepackage{dingbat}
\usepackage{wasysym}
\makeatletter
\g@addto@macro\appendix{%
  \addtocontents{toc}{\protect\setcounter{tocdepth}{1}}%
}
\makeatother

\usepackage{enumitem}

\begin{document}

\maketitle

\begin{abstract}
 
We prove a new concentration inequality for U-statistics of order two for uniformly ergodic Markov chains. Working with bounded and \textit{$\pi$-canonical} kernels, we show that we can recover the convergence rate of %\cite{arcones1993}
Arcones and Giné who proved a concentration result for U-statistics of independent random variables and canonical kernels. Our result allows for a dependence of the kernels~$h_{i,j}$ with the indexes in the sums, which prevents the use of standard blocking tools. Our proof relies on an inductive analysis where we use martingale techniques, uniform ergodicity, Nummelin splitting and Bernstein's type inequality. \\ Assuming further that the Markov chain starts from its invariant distribution, we prove a Bernstein-type concentration inequality that provides sharper convergence rate for small variance terms.
\end{abstract}

\section{Introduction}

Concentration of measure has been intensely studied during the last decades since it finds application in large span of topics such as model selection (see~\cite{massart07} and~\cite{Lerasle16}), statistical learning (see~\cite{clemencon17}), online learning (see~\cite{online-pairwise}) or random graphs (see~\cite{CL18} and~\cite{DCD20}). Important contributions in this field are those concerning U-statistics. A U-statistic of order~$m$ is a sum of the form \[\sum_{1\leq i_1<\dots<i_m\leq n} h_{i_1, \dots,i_m}(X_{i_1}, \dots,X_{i_m}),\] where~$X_1, \dots,X_n$ are independent random variables taking values in a measurable space~$(E,\Sigma)$ (with~$E$ Polish) and with respective laws~$P_i$ and where~$h_{i_1,\dots,i_m}$ are measurable functions of~$m$ variables~$h_{i_1,\dots,i_m}:E^{m} \rightarrow \mathds R$.

One important exponential inequality for U-statistics was provided by %Arcones and Giné in 
\cite{arcones1993} using a Rademacher chaos approach. Their result holds for bounded and canonical (or degenerate) kernels, namely satisfying for all~$ i_1,\dots,i_m\in [n]:=\{1,\dots,n\}$ with~$i_1<\dots< i_m$ and for all~$ x_1,\dots,x_m \in E$,
\[ 
\big\|h_{i_1, \dots, i_m}\big\|_{\infty}<\infty 
\quad
\text{and}
\quad
\forall j \in [1, n]\,,  
\ \mathds E_{X_j}
\Big[ h_{i_1,\dots,i_m}(x_{1}, \dots, x_{j-1},X_j,x_{j+1}, \dots ,x_m)
\Big]=0\,.\] They proved that in the degenerate case, the convergence rates for U-statistics are expected to be~$n^{m/2}$. Relying on precise moment inequalities of Rosenthal type, Giné, Latala and Zinn in~\cite{GineZinn} improved the result from~\cite{arcones1993}  by  providing  the  optimal four regimes of the tail, namely Gaussian, exponential, Weibull of orders 2/3 and 1/2. In the specific case of order 2 U-statistics, Houdré and Reynaud-Bouret in~\cite{houdre2002} recovered the result from~\cite{GineZinn} by replacing the moment estimates by martingales
type inequalities, giving as a by-product explicit constants. When the kernels are unbounded, it was shown that some results can be extended provided that the random variables~$h_{i_1, \dots,i_m}(X_{i_1},\dots,X_{i_m})$ have sufficiently light tails. One can mention~\cite[Theorem~3.26]{eichelsbacher03} where an exponential inequality for U-statistics with a single Banach-space valued, unbounded and canonical kernel is proved. Their approach is based on a decoupling argument originally obtained by~\cite{de1995decoupling} and the tail behavior of the summands is controlled by assuming that the kernel satisfies the so-called weak Cramér condition.   
It is now well-known that with heavy-tailed distribution for~$h_{i_1, \dots,i_m}(X_{i_1}, \dots,X_{i_m})$ we cannot expect to get exponential inequalities anymore. Nevertheless working with kernels that have finite~$p$-th moment for some~$p \in (1,2]$, Joly and Lugosi in~\cite{Lugosi15} construct an estimator of the mean of the U-process using the median-of-means technique that performs as well as the classical U-statistic with bounded kernels.\medskip

All the above mentioned results consider that the random variables~$(X_i)_{i\geq 1}$ are independent. This condition can be prohibitive for practical applications since modelization of real phenomena often involves some dependence structure. The simplest and the most widely used tool to incorporate such dependence is Markov chain. One can give the example of Reinforcement Learning (see~\cite{sutton11}) or Biology (see~\cite{suchard01}). Recent works provide extensions of the classical concentration results to the Markovian settings as~\cite{FJS18,jiang2018bernsteins,paulin15,adamczak:ustats,clemencon17}. The asymptotic behaviour of U-statistics in the Markovian setup has already been investigated by several papers. We refer to~\cite{bertail11} where the authors proved a Strong Law of Large Numbers and a Central
Limit Theorem proved for U-statistics of order 2 using the {\it renewal approach} based on the splitting technique. One can also mention~\cite{eichelsbacher01} regarding large deviation principles. However, there are only few results for the non-asymptotic behaviour of tails of U-statistics in a dependent framework.
The first results were provided in~\cite{Borisov2015} and~\cite{han16} where exponential inequalities for U-statistics of order~$m\geq2$ of time series under mixing conditions are proved. Those works were improved by~\cite{shen20} where a Hoeffding-type inequality for V and U-statistics is provided under conditions on the time dependent process that are easier to check in practice. In Section~\ref{apdx:han}, we describe in details the result of~\cite{shen20} and the differences with our work. Let us point out that all the above mentioned works regarding non-asymptotic tail bound for U-statistics in a dependent framework consider a fixed kernel, namely $h\equiv h_{i_1, \dots, i_m}$ for all $i_1,\dots,i_m$. Our work is the first to consider time dependent kernel functions which makes the theoretical analysis more challenging since the standard splitting method can be unworkable (cf. Section~\ref{sec:technical-assumption}). In Section~\ref{sec:appli-time-kernels} and~\ref{sec:time-dependent-cv-rate}, we stress the importance of working with index-dependent kernels for practical applications and we show on a specific example that one can reach significantly faster convergence rates with this approach. 

For the first time, we provide in this paper a Bernstein-type concentration inequality for U-statistics of order 2 in a dependent framework with kernels that may depend on the indexes of the sum and that are not assumed to be symmetric or smooth. We work on a general state space with bounded kernels that are~$\pi$-canonical. This latter notion was first introduced in~\cite{fort12} who proved a variance inequality for U-statistics of ergodic Markov chains. Our Bernstein bound holds for stationary chains but we provide a Hoeffding-type inequality without any assumption on the initial distribution of the Markov chain.

\subsection{Outline}

In Section~\ref{sec:assumptions}, we present and comment the assumptions under which our main results hold. In Section~\ref{subsec:csts}, we define and comment the key quantities involved in our results and we present our exponential inequalities with Theorems~\ref{mainthm} and~\ref{mainthm3} in Section~\ref{sec:expo-inequalities}. Section~\ref{sec:discussion} is dedicated to discussions where we give examples of Markov chains satisfying our assumptions and where we compare our results with the independent case. The proofs of both Theorems are presented in Section~\ref{sec:proofs}. In the Appendix, we provide the proof of some technical lemmas.

\clearpage
\section{Assumptions and notations}
\label{sec:assumptions}

We consider a Markov chain~$(X_i)_{i \geq 1}$ with transition kernel~$P:E\times E \rightarrow \mathds R$ taking values in a measurable space~$(E,\Sigma)$, and we introduce bounded functions~$h_{i,j}:E^2 \rightarrow \mathds R$. In this section, we describe the different assumptions on the Markov chain~$(X_i)_{i \geq 1}$ and on the functions $h_{i,j}$ that we will consider in Theorems~\ref{mainthm} and~\ref{mainthm3} presented in the next section.

\subsection{Uniform ergodicity}
\begin{assumption}
 \label{assumption1} The Markov chain~$(X_i)_{i\geq1}$ is~$\psi$-irreducible for some maximal irreducibility measure~$\psi$ on~$\Sigma$ (see~\cite[Section 4.2]{tweedie}). Moreover, there exist an integer~$m\geq1$,~$\delta_m>0$ and some probability measure~$\mu$ such that \[\forall x \in E, \; \forall A \in \Sigma, \quad \delta_m \mu(A) \leq P^m(x,A) .\]
 We denote by $\pi$ the unique invariant distribution of the Markov chain $(X_i)_{i\geq 1}$.
 \end{assumption}
 
For the reader familiar with the theory of Markov chains, Assumption~\ref{assumption1} states that the whole space~$E$ is a small set which is equivalent to the uniform ergodicity of the Markov chain~$(X_i)_{i \geq 1}$ (see~\cite[Theorem~16.0.2]{tweedie}), namely there exist constants~$0<\rho<1$ and~$L>0$ such that \[ \|P^n(x,\cdot)-\pi \|_{\mathrm{TV}} \leq L \rho^n, \qquad \forall n \geq0, \; \pi\mathrm{-a.e}\;  x \in E,  \]
where~$\pi$ is the unique invariant distribution of the chain~$(X_i)_{i\geq1}$ and where for any measure~$\omega$ on~$(E,\Sigma)$,~$\|\omega\|_{\mathrm{TV}}:= \sup_{A \in \Sigma}|\omega(A)|$ is the total variation norm of~$\omega$. From~\cite[section 2.3]{ferre12}), we also know that the Markov chain~$(X_i)_{i\geq 1}$ admits an absolute spectral gap~$1-\lambda >0$ with~$\lambda \in [0,1)$ (thanks to uniform ergodicity). We refer to~\cite[Section 3.1]{FJS18} for a reminder on the spectral gap of Markov chains.

\subsection{Upper-bounded Markov kernel} \label{sec:assumption2}

Assumption~\ref{assumption2} can be read as a reverse Doeblin's condition and allows us to achieve a change of measure in expectations in our proof to work with i.i.d. random variables with distribution~$\nu$. As a result, Assumption~\ref{assumption2} is the cornerstone of our approach since it allows to decouple the U-statistic in the proof.
\begin{assumption}
 \label{assumption2} There exist~$\delta_M>0$ and some probability measure~$\nu$ such that \[\forall x \in E,\; \forall A \in \Sigma, \quad  P(x,A) \leq \delta_M  \nu(A).\]  
 \end{assumption}

Assumption~\ref{assumption2} has already been used in the literature (see~\cite[Section 4.2]{lindsten2015}) and was introduced in~\cite{delmoral1999}. This condition can typically require the state space to be compact as highlighted in~\cite{lindsten2015}. 

Let us describe another situation where Assumption~\ref{assumption2} holds. Consider that~$(E,\|\cdot\|)$ is a normed space and that for all~$x \in E$,~$P(x,dy)$ has density~$p(x,\cdot)$ with respect to some measure~$\eta$ on~$(E,\Sigma)$. We further assume that there exists an integrable function~$u:E \to \mathds R_+$ such that $\forall x,y \in E, \quad p(x,y) \leq u(y) .$ Then considering for~$\nu$ the probability measure with density~$u/\|u\|_1$ with respect to~$\eta$ and~$\delta_M = \|u\|_1$,  Assumption~\ref{assumption2} holds.

\subsection{Exponential integrability of the regeneration time}
\label{splitting}

We introduce some additional notations which will be useful 
to apply Talagrand concentration result from~\cite{samson2000concentration}. Note that this section is inspired from~\cite{adamczak:ustats} and~\cite[Theorem~17.3.1]{tweedie}. We assume that Assumption~\ref{assumption1} is satisfied and we extend the Markov chain~$(X_i)_{i \geq 1}$ to a new (so called \textit{split}) chain~$(\widetilde{X}_n, R_n)\in E\times \{0,1\}$ (see~\cite[Section 5.1]{tweedie} for a reminder on the splitting technique), satisfying the following properties.

\begin{itemize}
\item~$(\widetilde{X}_n)_n$ is again a Markov chain with transition kernel~$P$ with the same initial distribution as~$(X_n)_n$. We recall that~$\pi$ is the invariant distribution on the~$E$.
\item if we define~$T_1=\inf \{n>0 : R_{nm}=1\}$, \[T_{i+1}= \inf \{n>0 : R_{(T_1+\dots+T_i+n)m}=1\},\]
then~$T_1, T_2,\dots$ are well defined and independent. Moreover~$T_2,T_3,\dots$ are i.i.d.
\item if we define~$S_i=T_1+\dots+T_i$, then the ‘‘blocks'' \[Y_0= (\widetilde X_1,\dots ,\widetilde X_{mT_1+m-1}),\quad \text{and}\quad Y_i= (\widetilde X_{m(S_i+1)},\dots ,\widetilde X_{m(S_{i+1}+1)-1}), \  i >0,\]form a one-dependent sequence (i.e. for all~$i$,~$\sigma((Y_j)_{j<i})$ and~$\sigma((Y_j)_{j>i})$ are independent). Moreover, the sequence~$Y_1, Y_2, \dots$ is stationary and if~$m=1$ the variables ~$Y_0, Y_1,\dots$ are independent. In consequence, for any measurable space~$(S,\mathcal B)$ and measurable functions~$f:S \rightarrow \mathds R$, the variables \[Z_i=Z_i(f) =\sum_{j=m(S_i+1)}^{m(S_{i+1}+1)-1}f(\widetilde X_j), \quad i\geq1,\]constitute a one-dependent sequence (an i.i.d.  sequence if~$m= 1$). Additionally, if~$f$ is~$\pi$-integrable (recall that~$\pi$ is the unique stationary measure for the chain), then\[ \mathds E[Z_i]=\delta_m^{-1}m\int f d\pi.\]
\item the distribution of~$T_1$ depends only on~$\pi$,~$P$,~$\delta_m$,~$\mu$, whereas the law of~$T_2$ only on~$P$,~$\delta_m$ and~$\mu$.
\end{itemize}

{\bf Remark} Let us highlight that~$(\widetilde{X}_n)_n$ is a Markov chain with transition kernel~$P$ and same initial distribution as~$(X_n)_n$. Hence for  our  purposes  of  estimating  the  tail  probabilities, we will identify~$(X_n)_n$ and~$(\widetilde{X}_n)_n$.

To derive a concentration inequality, we use the exponential integrability of the regeneration times which is ensured if the chain is uniformly ergodic as stated by Proposition~\ref{prop:tau}. A proof can be found in Section \ref{proof:prop-regeneration} of the Appendix.

\begin{definition}\label{orlicz}
For~$\alpha>0$, define the function~$\psi_{\alpha}:\mathds R_+\rightarrow \mathds R_+$ with the formula~$\psi_{\alpha}(x) =\exp(x\alpha)-1$.  Then for a random variable~$X$, the~$\alpha$-Orlicz norm is given by \[\|X\|_{\psi_{\alpha}}= \inf\left\{\gamma >0 :\mathds E [\psi_{\alpha}(|X|/\gamma)]\leq1\right\}.\]
\end{definition}

\begin{proposition} \label{prop:tau} If Assumption~\ref{assumption1} holds, then
 \begin{equation}\|T_1\|_{\psi_1}<\infty \quad \text{and}\quad \|T_2\|_{\psi_1}<\infty, \label{eq:regeneration}\end{equation}
 where~$\|\cdot\|_{\psi_1}$ is the~$1$-Orlicz norm introduced in Definition~\ref{orlicz}. We denote~$\tau:=\max(\|T_1\|_{\psi_1},\|T_2\|_{\psi_1}).$
\end{proposition}

\subsection{$\pi$-canonical and bounded kernels}
With Assumption~\ref{assumption4}, we introduce the notion of \textit{$\pi$-canonical} kernel which is the counterpart of the canonical property from~\cite{ginenick}.

\begin{assumption} \label{assumption4} Let us denote~$\mathcal B(\mathds R)$ the Borel algebra on~$\mathds R$. For all~$i,j \in [n]$, we assume that~$h_{i,j}:(E^2,\Sigma \otimes \Sigma) \to (\mathds R,\mathcal B(\mathds R))$ is measurable and is~$\pi$-canonical, namely \[\forall x,y \in E, \quad \mathds E_{\pi}[h_{i,j}(X,x)]=\mathds E_{\pi}[h_{i,j}(X,y)] = \mathds E_{\pi}[h_{i,j}(x,X)]=\mathds E_{\pi}[h_{i,j}(y,X)].\]
 This common expectation will be denoted~$\mathds E_{\pi}[h_{i,j}]$. \\ Moreover, we assume that for all~$i,j \in [n]$,~$\|h_{i,j}\|_{\infty} <\infty.$
\end{assumption}

\noindent
\textbf{Remarks}

\noindent
$\bullet$ A large span of kernels are~$\pi$-canonical. This is the case of translation-invariant kernels which have been widely studied in the Machine Learning community. Another example of~$\pi$-canonical kernel is a rotation invariant kernel when~$E=\mathds S^{d-1}:=\{x \in \mathds R^d : \|x\|_2=1\}$ with~$\pi$ also rotation invariant (see~\cite{CL18} or~\cite{DCD20}).  

\smallskip

\noindent
$\bullet$ The notion of~$\pi$-canonical kernels is the counterpart of canonical kernels in the i.i.d. framework (see for example~\cite{houdre2002}). Note that we are not the first to introduce the notion of~$\pi$-canonical kernels working with Markov chains. In~\cite{fort12}, Fort and al. provide a variance inequality for U-statistics whose underlying sequence of random variables is an ergodic Markov Chain. Their results holds for~$\pi$-canonical kernels as stated with~\cite[Assumption A2]{fort12}.

\smallskip

\noindent
$\bullet$ Note that if the kernels~$h_{i,j}$ are not~$\pi$-canonical, the U-statistic decomposes into a linear term and a~$\pi$-canonical U-statistic. This is called the \textit{Hoeffding decomposition} (see~\cite[p.176]{ginenick}) and takes the following form 
\begin{align*}
&\sum_{i\neq j}\left( h_{i,j}(X_i,X_j)- \mathds{E}_{(X,Y) \sim \pi \otimes\pi}[h_{i,j}(X,Y)]\right)=\sum_{i\neq j}\widetilde h_{i,j}(X_i,X_j) - \mathds E_{ \pi } \left[\widetilde h_{i,j} \right] \\
&\qquad \qquad + \sum_{i\neq j} \left(\mathds E_{X \sim \pi} \left[ h_{i,j}(X,X_j) \right] - \mathds E_{(X,Y) \sim \pi \otimes\pi} \left[ h_{i,j}(X,Y) \right]\right)\\
&\qquad \qquad +  \sum_{i\neq j} \left(\mathds E_{X \sim \pi} \left[ h_{i,j}(X_i,X) \right] - \mathds E_{(X,Y) \sim \pi \otimes\pi} \left[ h_{i,j}(X,Y) \right]\right),
\end{align*}
where for all~$ j$, the kernel~$\widetilde h_{i,j}$ is~$\pi$-canonical with 
\[\forall x,y \in E, \quad\widetilde h_{i,j}(x,y) = h_{i,j}(x,y)- \mathds E_{X \sim \pi} \left[ h_{i,j}(x,X) \right]-\mathds E_{X \sim \pi} \left[ h_{i,j}(X,y) \right].\]

\subsection{Additional technical assumption}
\label{sec:technical-assumption}

In the case where the kernels~$h_{i,j}$ depend on both~$i$ and~$j$, we need Assumption~\ref{assumption0}$.(ii)$ to prove Theorem~\ref{mainthm}. Assumption~\ref{assumption0}.$(ii)$ is a mild condition on the initial distribution of the Markov chain that is used when we apply Bernstein's inequality for Markov chains from Proposition~\ref{bernstein-markov} (see Section \ref{apdx:EZ-chi} of the Appendix).

\begin{assumption}
\label{assumption0} At least one of the following conditions holds.

$(i) \;$~For all~$i,j \in [n]$,~$h_{i,j}\equiv h_{1,j},$ i.e. the kernel function~$h_{i,j}$ does not depend on~$i.$

$(ii)\;~$The initial distribution of the Markov chain~$(X_i)_{i\geq1}$, denoted~$\chi$, is absolutely continuous with respect to the invariant measure~$\pi$ and its density~$\frac{d\chi}{d\pi}$ has finite~$p$-moment for some~$p\in(1,\infty]$, i.e
\[\infty > \left\|\frac{d\chi}{d\pi}\right\|_{\pi,p} := \left\{
    \begin{array}{ll}
      \left[ \int \left| \frac{d\chi}{d\pi} \right|^p d\pi \right]^{1/p} & \mbox{if } p<\infty, \\
    \mathrm{ess }\sup \; \left| \frac{d\chi}{d\pi}\right|& \mbox{if }p=\infty.
    \end{array}
\right.\]
In the following, we will denote~$q = \frac{p}{p-1}\in [1, \infty)$ (with~$q=1$ if~$p=+\infty$) which satisfies~$\frac{1}{p}+\frac{1}{q}=1$.
\end{assumption}

Assumption~\ref{assumption0} is needed at one specific step of our proof where we need to bound with high probability \[ \sum_{j=2}^n\mathds E\big[\big| \sum_{i=1}^{j-1}p_{i,j}(X_i,X'_j)\big|^k\big],\quad \text{with} \quad \forall i,j ,\quad \forall x,y\in E,\quad p_{i,j}(x,y) := h_{i,j}(x,y)-\mathds E_{\pi} [h_{i,j}] ,\]
and where $(X'_j)_{j\geq 1}$ are i.i.d. random variables with distribution $\nu$ from Assumption~\ref{assumption2}. In the case where Assumption~\ref{assumption0}$.(i)$ holds, we can use for any fixed $j \in \{2,\dots,n\}$ the splitting method to decompose the sum $\sum_{i=1}^{j-1}p_{i,j}(X_i,X'_j)$ in different blocks whose lengths are given by the regeneration times of the split chain. Thanks to Assumption~\ref{assumption0}$.(i)$, those blocks are independent and we can use standard concentration tools for sums of independent random variables. This approach is valid for any initial distribution of the chain. However, if Assumption~\ref{assumption0}$.(i)$ is not satisfied, the blocks used to decompose $\sum_{i=1}^{j-1}p_{i,j}(X_i,X'_j)$ are not independent and the splitting method can no longer be used. To bypass this issue, we need a Bernstein-type concentration inequality for additive functionals of Markov chains with  time-dependent functions (see Proposition~\ref{bernstein-markov} in Section \ref{apdx:EZ-chi} of the Appendix). %~\ref{bernstein-markov} in Section~\ref{apdx:EZ-chi}). Proposition~\ref{bernstein-markov} 
Proposition~4 is a straightforward corollary of~\cite[Theorem 1]{jiang2018bernsteins} and requires Assumption~\ref{assumption0}.$(ii)$ to be satisfied. We refer to Section~\ref{proof:mainproposition} and in particular to Section~\ref{proof:EZK} for further details.

\clearpage

\section{Main results}
\label{sec:main-results}

\subsection{Preliminary comments}
\label{subsec:csts}

Under the assumptions presented in Section~\ref{sec:assumptions}, Theorem~\ref{mainthm} and~\ref{mainthm3} provided in Section~\ref{sec:expo-inequalities} give exponential inequalities for the U-statistic \[ U_{\mathrm{stat}}(n)=\sum_{1\leq i<j \leq n}\left( h_{i,j}(X_i,X_j)- \mathds E\left[h_{i,j}( X_i, X_j)\right]\right).\]
Theorem~\ref{mainthm} provides an Hoeffding-type concentration result that holds without any (or mild) condition on the initial distribution of the chain. By assuming that the chain $(X_i)_{i\geq1}$ is stationary (meaning that $X_1$ is distributed according to $\pi$), Theorem~\ref{mainthm3} gives a Bernstein-type concentration inequality and leads to a better convergence rate compared to Theorem~\ref{mainthm}.

The proof of our main results relies on a martingale technique conducted by induction at depth 
$t_n := \floor{r\log n}$ with~$r > 2\left(\log(1/\rho)\right)^{-1}$ (see the remark following Assumption~\ref{assumption1} for the definition of $\rho$). With the notations of Section~\ref{sec:assumptions}, our concentration inequalities involve the following quantities
\begin{align}A:=& 2\max_{i,j} \|h_{i,j}\|_{\infty}, 
\quad C_n^2:=\sum_{j=2}^n \sum_{i=1}^{j-1} \mathds E \left[\mathds E_{X' \sim \nu} [ p_{i,j}^2(X_i,X')] \right], \\
B_n^2:=& \max\Bigg[ \max_{0 \leq k\leq t_n} \; \max_{i} \;\sup_x \sum_{j=i+1}^n \mathds E_{X' \sim \nu} \left(\mathds E_{X \sim P^k(X',\cdot)} \; p_{i,j}(x,X) \right)^2 , \notag\\
&\qquad  \quad \max_{0 \leq k\leq t_n}\; \max_{j}\; \sup_{y}\sum_{i=1}^{j-1}\mathds E_{\tilde X \sim \pi}\left( \mathds E_{X \sim P^k(y,\cdot)} p_{i,j}(\tilde X, X)\right)^2 \Bigg],\label{eq:defB}\end{align}
with the convention that~$P^0(y,\cdot)$ is the Dirac measure at point~$y\in E$. Let us comment those terms.

\begin{itemize}
\item \underline{Understanding of the origin of $B_n$.} $B_n$ involves supremums over $k$ ranging from $0$ to $t_n$. The terms in the supremum corresponding to some specific value of $k$ arise in our proof at the $k$-th step of our induction procedure (and will be denoted $\mathfrak B_k$ in Section~\ref{sec:proofs}, so that $B_n=\sup_{0\leq k\leq t_n }\mathfrak B_k$).
\item \underline{Bounding $B_n$ with uniform ergodicity.} The uniform ergodicity of the Markov chain ensured by Assumption~\ref{assumption1} can allow to bound~$B_n$ since for all~$x,y \in E$ and for all~$k\geq 0$, \[ \left|  \mathds E_{X \sim P^k(y,\cdot)} p_{i,j}(x,X)\right| \leq \sup_z |h_{i,j}(x,z)|\times  \|P^k(y,\cdot)-\pi\|_{\mathrm{TV}} .\]

\item \underline{The case where $\nu=\pi$ and the independent setting}\\
In the specific case where~$\nu=\pi$ (which includes the independent setting), we get that
\[C_n^2= \sum_{i<j} \mathds E\left\{ \mathrm{Var}_{\tilde X\sim \pi}\left[ h_{i,j}(X_i,\tilde X)|X_i\right] \right\},\]
and using Jensen inequality that
\[ B_n^2\leq \max\Big[\sup_{x,i}\sum_{j=i+1}^n \mathrm{Var}_{\tilde X\sim \pi}[h_{i,j}(x,\tilde X)] \; ,\;  \sup_{y,j} \; \sum_{i=1}^{j-1} \mathrm{Var}_{\tilde X\sim \pi}[h_{i,j}(\tilde X,y)] \Big].\]
Hence,~$C_n^2$ and~$B_n^2$ can be understood as variance terms that would tend to be larger as~$\nu$ moves away from~$\pi$. Let us point out that in the independent setting, all terms for $k$ ranging from 1 to $t_n$ in the definition of $B_n^2$ vanish but the term corresponding to $k=0$ does not since $P^0(y,\cdot)$ is the Dirac measure at $y$. We provide a detailed comparison of our results with known exponential inequalities in the independent setting in Section~\ref{sec:compa-inde}.
\item \underline{Bounding~$B_n$ and~$C_n$: A way to read immediately the convergence rates in our main results}\\
Using coarse bounds, one immediately gets that~$B_n \leq A\sqrt n$ and~$C_n \leq An$. We prompt the reader to keep in mind these bounds in order to directly see the rate of convergence and the dominant terms in the inequalities from Section~\ref{sec:expo-inequalities}. These bounds can be significantly improved for particular cases as done in the example presented in Section~\ref{sec:time-dependent-cv-rate}.
\end{itemize}

\subsection{Exponential inequalities}
\label{sec:expo-inequalities}

We now state our first result Theorem~\ref{mainthm} whose proof can be found in Section~\ref{sec:proof-mainthm}.

\begin{theorem}
\label{mainthm} Let~$n \geq 2$. We suppose Assumptions~\ref{assumption1},~\ref{assumption2} and~\ref{assumption4} described in Section~\ref{sec:assumptions}. 
There exist two constants~$\beta,\kappa>0$ such that for any~$u>0$, 

\begin{enumerate}[label=\alph*)]
\item if Assumption~\ref{assumption0}$.(i)$ is satisfied, it holds with probability at least~$1-\beta e^{-u}\log(n),$
\[\hspace{-0.3cm} U_{\mathrm{stat}}(n)\leq  \kappa \log(n)\Big( \left[A \log(n) \sqrt n\right]\sqrt u + \left[A+  B_n\sqrt n\right] u +\left[2A\sqrt n \right]u^{3/2} + A \left[u^2+n\right]  \Big),\]
\item if Assumption~\ref{assumption0}$.(ii)$ is satisfied, it holds with probability at least~$1-\beta e^{-u}\log(n),$
\[\hspace{-0.3cm} U_{\mathrm{stat}}(n)\leq  \kappa \log(n)\Big( \left[ \mathbf C_n+A \log(n)\sqrt{n} \right]\sqrt u + \left[A+  B_n\sqrt n\right] u +\left[2A\sqrt n \right]u^{3/2} + A \left[u^2+n\right]  \Big).\]
\end{enumerate}

\end{theorem}

Note that the kernels~$h_{i,j}$ do not need to be symmetric and that we do not consider any assumption on the initial measure of the Markov chain~$(X_i)_{i\geq1}$ if the kernels~$h_{i,j}$ do not depend on~$i$ (see Assumption~\ref{assumption0}). %Let us mention that depending on whether Assumption~\ref{assumption0}.$(i)$ or Assumption~\ref{assumption0}.$(ii)$ holds, we do not obtain exactly the same bounds and we take the maximum of both results in Theorem~\ref{mainthm}.
By bounding coarsely~$B_n$ and~$C_n$ in Theorem~\ref{mainthm} (respectively by $\sqrt n A$ and $nA$), we get that there exist constants~$\beta,\kappa>0$ such that for any~$u \geq 1$, it holds with probability at least~$1-\beta e^{-u}\log n$,
\begin{equation}\label{eq:mainthm2}\frac{2}{n(n-1)} U_{\mathrm{stat}}(n) \leq \kappa \max_{i,j}\|h_{i,j}\|_{\infty}\log n\,\, \bigg\{   \frac{u}{n} + \left[\frac{u}{n}\right]^{2}  \bigg\}. \end{equation}
In particular it holds
\[
\frac{2}{n(n-1)} U_{\mathrm{stat}}(n)
=\mathcal O_{\mathds P}
\Bigg(\frac{\log (n)\log\log n}{n}\Bigg)\,,
\]
where~$\mathcal O_{\mathds P}$ denotes stochastic boundedness. Up to a~$\log (n)\log\log n$ multiplicative term, we uncover the optimal rate of Hoeffding’s inequality for canonical U-statistics of order~$2$, see~\cite{Lugosi15}. Taking a close look at the proof of Theorem~\ref{mainthm} (and more specifically at Section~\ref{proof:mainproposition-remaining}), one can remark that the same results hold if the U-statistic is centered with the expectations~$\mathds E_{\pi}[h_{i,j}]$, namely for \[\sum_{1\leq i<j\leq n} \left(  h_{i,j}(X_i,X_j)-\mathds E_{\pi}[h_{i,j}]\right).\] 
It is well-known that one can expect a better convergence rate when variance terms are small with a Bernstein bound. The main limitation in Theorem~\ref{mainthm} that prevents us from taking advantage of small variances is the term at the extreme right on the concentration inequality of Theorem~\ref{mainthm}, namely~$An \log n$. Working with the additional assumption that the Markov chain~$(X_i)_{i\geq 1}$ is stationary -- meaning that the initial distribution of the chain is the invariant distribution~$\pi$ -- we are able to prove a Bernstein-type concentration inequality as stated with Theorem~\ref{mainthm3}. The proof of Theorem~\ref{mainthm3} is provided in Section~\ref{proof-mainthm3}. Stationarity is only used to bound the remaining terms that were not already considered in the $t_n$ steps of our induction procedure (see Section~\ref{subsec:csts} for the definition of $t_n$). We refer to the proof of Proposition~\ref{mainproposition-remaining}$.b)$ in Section~\ref{proof:mainproposition-remaining} for details.

\begin{theorem} \label{mainthm3}  We suppose Assumptions~\ref{assumption1},~\ref{assumption2} and~\ref{assumption4} described in Section~\ref{sec:assumptions}. We further assume that the Markov chain~$(X_i)_{i\geq 1}$ is stationary. Then there exist two constants~$\beta,\kappa>0$ such that for any~$u>0$, it holds with probability at least~$1-\beta e^{-u}\log n$,
\[U_{\mathrm{stat}}(n)\leq  \kappa \log(n)\Big( \left[ C_n+A \log(n)\sqrt{n} \right]\sqrt u + \left[A+  B_n\sqrt n\right] u +\left[2A\sqrt n \right]u^{3/2} + A \left[u^2+\mathbf{\log n}\right]  \Big).\]
In case where Assumption~\ref{assumption0}$.(i)$ holds, one can remove~$C_n$ in the previous inequality.
\end{theorem}

\clearpage

\subsection{Discussion}

\label{sec:discussion}

\subsubsection{Examples of Markov chains satisfying the Assumptions}
\label{subsec:examples}

\paragraph{Example 1: Finite state space.}
For Markov chains with finite state space, Assumption~\ref{assumption2} holds trivially. Hence, in such framework the result of Theorem~\ref{mainthm} holds for any uniformly ergodic Markov chain. In particular, this is true for any aperiodic and irreducible Markov chains using~\cite[Lemma 7.3.(ii)]{behrends2000introduction}.

%Assume that the space~$E$ is compact. If the Markov kernel is absolutely continuous with respect to a given probability measure~$\nu$ on~$E$ with~$P(x,dy) = f(x,y)\nu(dy)$ and~$\sup f(x,\cdot)$ continuous for all~$x \in E$~\ref{assumption2}.(i) are trivially satisfy  and  

\paragraph{Example 2: AR(1) process.}
Let us consider the process~$(X_n)_{n \in \mathds N}$ on~$\mathds R^k$ defined by \[ X_0 \in \mathds R^k \text{  and for all  } n  \in \mathds N, \quad X_{n+1}=H(X_n)+Z_n,\] where~$(Z_n)_{n \in \mathds N}$ are i.i.d random variables in~$\mathds R^k$ and~$H: \mathds R ^k \to \mathds R^k$ is an application. Such a process is called an auto-regressive process of order 1, noted AR(1). Assuming that the distribution of $Z_1$ has density $f_Z$ with respect to the Lebesgue measure on $\mathds R^k$ (denoted $\lambda_{Leb}$), it is well-known that mild regularity conditions on $H$ and $f_Z$ ensure that the Markov chain $(X_i)_{i\geq1}$ is uniformly ergodic. These conditions require in particular that both $f_Z$ and $H$ are continuously differentiable with $H$ bounded. We refer to~\cite{Doukhan-Ghindes} for the full statement.

We denote $B_H:=B(0,\|H\|_{\infty})$ the euclidean ball in $\mathds R^k$ with radius $\|H\|_{\infty}$ centered at $0$. Assuming that~$y \mapsto {\sup}_{ \{z \in B_H\}} f_Z(y-z)$ is integrable on~$\mathds R^k$ with respect to~$\lambda_{Leb}$, we get that  Assumption~\ref{assumption2} holds (see the remark following Assumption~\ref{assumption2}). The previous condition on~$f_Z$ is for example satisfied for Gaussian distributions. We deduce that Theorems~\ref{mainthm} and~\ref{mainthm3} can be applied in such settings that are typically found in nonlinear filtering problem (see~\cite[Section 4]{delmoral1999}).

 \paragraph{Example 3: ARCH process.}
 
 Let us consider~$E=\mathds R$. The ARCH model is \[X_{n+1} = H(X_n)+G(X_n)Z_{n+1},\]
 where~$H$ and~$G$ are continuous functions, and~$(Z_n)_n$ are i.i.d. centered normal random variables with variance~$\sigma^2>0$. Assuming that~$\inf_x |G(x)|\geq a>0$, we know that the Markov chain~$(X_n)_n$ is irreducible and aperiodic (see~\cite[Lemma 1]{ANGONZE98}). Assuming further that~$\|H\|_{\infty}\leq b < \infty$ and that~$\|G\|_{\infty}\leq c$, we can show that Assumptions~\ref{assumption1} and~\ref{assumption2} hold. Let us first remark that the transition kernel~$P$ of the Markov chain~$(X_n)_n$ is such that for any~$x \in \mathds R$,~$P(x, dy)$ has density~$p(x,\cdot)$ with respect to the Lebesgue measure with \[p(x,y)=\left( 2 \pi \sigma^2\right)^{-1} \exp \left( - \frac{(y-H(x))^2}{2\sigma^2 G(x)^2} \right).\] 
Defining for any $y\in \mathds R$,
 \[g_m(y):= \frac{1}{2 \pi \sigma^2} \times \left\{
    \begin{array}{ll}
      \exp \left( - \frac{(y-B)^2}{2\sigma^2 a^2} \right) & \mbox{if } y< -b \\
     \exp \left( - \frac{2B^2}{\sigma^2 a^2} \right) & \mbox{if } |y| \leq b\\
      \exp \left( - \frac{(y+B)^2}{2\sigma^2 a^2} \right) & \mbox{if } y>b     \end{array}
\right. \quad \text{and} \quad  g_M(y) :=\left( 2 \pi \sigma^2\right)^{-1} \times \left\{
    \begin{array}{ll}
      \exp \left( - \frac{(y+b)^2}{2\sigma^2 c^2} \right) & \mbox{if } y< -b \\
    1 & \mbox{if } |y| \leq b\\
      \exp \left( - \frac{(y-b)^2}{2\sigma^2 c^2} \right) & \mbox{if } y>b   \end{array}
\right.,\]
it holds $g_m(y)\leq p(x,y)\leq g_M(y)$ for any $x,y\in \mathds R$. We deduce that considering~$\delta_m=\|g_m\|_1, \delta_M=\|g_M\|_1$ and~$\mu$ (resp.~$\nu$) with density~$g_m/\|g_m\|_1$ (resp.~$g_M/\|g_M\|_1$) with respect to the Lebesgue measure on~$\mathds R$, Assumptions~\ref{assumption1} and~\ref{assumption2} hold.

\subsubsection{The independent setting}
\label{sec:compa-inde}

In the case where the random variables $(X_i)_{i\geq 1}$ are independent, the terms~$B_n^2$ and~$C_n^2$ involved in our exponential inequality take the following form 
\begin{align}\label{eq:BCinde}
 B_n^2&=\max \left[ \sup_{i,x} \sum_{j=i+1}^n \mathds E\left[ p_{i,j}^2(x,X_j) \right]\; , \; \sup_{j,y} \sum_{i=1}^{j-1} \mathds E\left[ p_{i,j}^2(X_i,y) \right] \right]\\
\quad \text{and}\quad C_n^2&= \sum_{j=2}^n \sum_{i=1}^{j-1} \mathds E\left[ p_{i,j}(X_i,X_j)^2\right],\notag
\end{align}
where we remind that all terms for $k\in \{1,\dots,t_n\}$ in the definition of $B_n^2$ in Eq.\eqref{eq:defB} vanish and it only remains the contribution of terms for $k=0$. 
In the independent setting, \cite[Theorem 1]{houdre2002} proved that for any $u>0$, it holds with probability at least $1-3e^{-u}$,
\[U_{\mathrm{stat}}(n) \leq C_n\sqrt u + \left(D_n+F_n\right)u+B_nu^{3/2}+Au^2,\]
where $A$, $B_n$ and $C_n$ coincide with the terms of this paper (see Eq.\eqref{eq:BCinde}). Let us comment the quantities involved in the different regimes of the tail behaviour.
\begin{itemize}
\item \underline{Sub-Gaussian.} In Theorem~\ref{mainthm3}, we recover the term~$C_n$ from~\cite{houdre2002} and we suffer an additional $A\sqrt n \log n$ term.
\item \underline{Sub-Exponential.} $D_n$ and~$F_n$ come from duality arguments in the proof of~\cite{houdre2002}. We do not recover the counterpart of these terms in Theorem~\ref{mainthm3} since working with dependent variables bring additional technical difficulties and the use for example of a decoupling argument. $D_n+F_n$ is replaced by $A+B_n\sqrt n$ in our result.
\item \underline{Sub-Weibull with parameter $2/3$.} While~\cite{houdre2002} finds the quantity~$B_n$ for the term $u^{3/2}$, the counterpart in Theorem~\ref{mainthm3} is the worst case scenario since it always holds $B_n \leq A \sqrt n$.
\item \underline{Sub-Weibull with parameter $1/2$.} We obtain the same behaviour for the sub-Weibull (with parameter $1/2$) regime of the tail behaviour.
\end{itemize}
Let us also mention that Theorem~\ref{mainthm3} has an additive term $A\log^2 n$ (that will not be dominant for standard choice of $u$). This term can be understood as a proof artefact and arises when we bound the remaining terms in the U-statistic that were not considered in our induction procedure. We finally point out that our result involves additive $\log n$ factors (both in the tail bound and in the probability).

\subsection{Connections with the literature}

\label{apdx:han}

In this section, we describe the concentration inequality obtained in~\cite{shen20} for U-statistics in a dependent framework and we explain the differences with our work. We consider an integer~$n \in \mathds N \backslash \{0\}$ and a geometrically~$\alpha$-mixing sequence $(X_i)_{i\in [n]}$ (see~\cite[Section 2]{shen20}) with coefficient \[\alpha(i)\leq \gamma_1 \exp(-\gamma_2i), \quad \text{for all }i\geq 1,\] where~$\gamma_1,\gamma_2$ are two positive absolute constants. We consider a kernel~$h:\mathds R^d \times \mathds R^d \to \mathds R$ degenerate, symmetric, continuous, integrable and satisfying for some~$q\geq 1,$~$\int _{\mathds R^{2d}} |\mathcal F h(u)| \|u\|^q_2 du<\infty$, where $\mathcal F h$ denotes the Fourier-transform of $h$. Then Eq.(2.4) from~\cite{shen20} states that there exists a constant~$c>0$ such that for any~$u >0$, it holds with probability at least~$1-6e^{-u}$
\begin{equation} \label{shen20:ustat}\frac{2}{n(n-1)} U_{\mathrm{stat}}(n) \leq 4c\|\mathcal F h\|_{L^1} \left\{ A_n^{1/2} \frac{u}{n} +  c\log^4(n) \left[\frac{u}{n}\right]^{2}\right\},\end{equation}
where~$A_{n}^{1/2} = 4  \left( \frac{64 \gamma_1^{1/3}}{1-\exp(-\gamma_2/3)} + \frac{\log^4(n)}{n}  \right)$ and~$ U_{\mathrm{stat}}(n)=\sum_{1\leq i<j \leq n}\left( h(X_i,X_j)- \mathds{E}_{\pi}[h]\right)$.

\cite{shen20} has the merit of working with geometrically~$\alpha$-mixing stationary sequences which includes in particular geometrically (and hence uniformly) ergodic Markov chains (see~\cite[p.6]{jones2004}). For the sake of simplicity, we presented the result of~\cite{shen20} for U-statistics of order 2, but their result holds for U-statistics of arbitrary order~$m \geq 2$. Nevertheless, they only consider state spaces like~$\mathds R^d$ with~$d \geq 1$ and they work with a unique kernel~$h$ (i.e.~$h_{i_1, \dots ,i_m}=h$ for any~$i_1, \dots,i_m$) which is assumed to be symmetric continuous, integrable and that satisfies some smoothness assumption. On the contrary, we consider general state spaces and we allow different kernels $h_{i,j}$ that are not assumed to be symmetric or smooth. In addition, Theorem~\ref{mainthm3} is a Bernstein-type exponential inequality where we can benefit from small variance terms, which is not the case for~\cite{shen20}. We provide a specific example in Section~\ref{sec:time-dependent-cv-rate}.

\subsubsection{Motivations for the study of time dependent kernels}
\label{sec:appli-time-kernels}

In this section, we want to stress the importance of working with weighted U-statistics for practical applications. In the following, we detail two specific examples borrowed from the fields of information retrieval and of homogeneity tests. Note that one could find other applications such as in genetic association (cf.~\cite{geneteicweightedUstat}) or for independence tests (cf.~\cite{shieh}). 
\paragraph{Average-Precision Correlation}
When we search the Internet, the browser computes a numeric score on how well each object in the database matches the query, and rank the objects according to this value. In order to evaluate the quality of this browser, a standard approach in the field of information retrieval consists in comparing the ranking provided by the web search engine and the ranking obtained from human labels (cf.~\cite{inforetrievalUstat}). One way to measure how well both rankings are aligned is to report the correlation between them. One of the most commonly used rank correlation statistic is the Kendall’s~$\tau$. Considering a dataset of size~$n \in \mathds N$ ordered according to the human labels and denoting~$X_i$ the rank the browser gives to the~$i$-th element, the Kendall's~$\tau$ is defined by
\[\tau^{\mathrm{Ken}}:=\frac{2}{n(n-1)} \sum_{i\neq j} \left\{\mathds 1_{X_i > X_j } \mathds 1_{i>j} + \mathds 1_{X_i<X_j} \mathds 1_{i<j} \right\}- 1. \]
Since only the top ranking objects are shown to the user, it would be legitimate to penalize heavier errors made on items having high rankings. The Kendall's~$\tau$ does not make such distinctions and new correlation measurements have been popularized to address this issue. One of them is the so-called Average-Precision Correlation~(cf.~\cite{inforetrieval}) which is defined by
\[\tau^{\mathrm{AP}}:= \frac{2}{n-1} \sum_{j=2}^n \frac{\sum_{i=1}^{j-1} \mathds 1_{X_i < X_j}}{j-1}
- 1. \]
Note that $\tau^{\mathrm{AP}}$ is a U-statistic where the kernels~$h_{i,j}(x,y): = \frac{\mathds 1_{x < y}}{j-1}$ depend on~$j$. Let us point out that~$h_{i,j}$ do not depend on~$i$ so that Assumption~\ref{assumption0}$.(i)$ holds.

\paragraph{Accounting for confounding covariates}
U-statistics are powerful tools to compare the distributions of random variables across two groups (say with labels~$0$ and~$1$) from samples~$X_1,\dots,X_n$ and~$X_{n+1},\dots,X_{n+m}$. The typical example is the Wilcoxon Rank Sum Test (WRST) based on the following U-statistic
\[\sum_{i=1}^n\sum_{j=1}^m h(X_i,X_{n+j}) \quad \text{where} \quad h(x,y):= \frac12\mathds 1_{x<y}+\frac12 \mathds 1_{x\leq y}.\]
The WRST relies on the following idea: if the data is pooled and then ranked, the
average rank of observations from each group should be the same. For any~$i\in[n+m]$, let~$G_i$ be the random variable valued in~$\{0,1\}$ allocating the~$i$-th individual to one of the two groups. Note that the observed allocation~$(g_i)_{i\in[n+m]}$ is given by~$g_i=0$ if and only if~$i\leq n$. When group membership is not assigned through randomization, there may be confounding covariates~$Z$ (assumed to be observed) that can cause a spurious association between outcome and group membership. In that case, we wish rather to test the null hypothesis~$\mathds P(X \leq t \, |\, G=0, Z = z) = \mathds P(X \leq t \, |\, Z=z)$.  In~\cite{confoundingUstat}, the authors developed such a test by working with the following adjusted U-statistics involving index-dependent kernels
\[\big(\sum_{i\,:\,g_i=0}w(z_i,g_i) \sum_{j\,:\,g_j=1}w(z_j,g_j) \big)^{-1}\sum_{i\,:\,g_i=0} \sum_{j\;:\,g_j=1} h(X_i,X_j) w(z_i,g_i)w(z_j,g_j), \]
where the weights $w(z_i, g_i) = (\mathds P(G=g_i\, |\, Z=z_i))^{-1}$ can be estimated with a logistic regression.

\subsubsection{Time dependent kernels and convergence rate}
\label{sec:time-dependent-cv-rate}

In this section, we consider a stationary Markov chain $(X_i)_{i\geq1}$ satisfying Assumptions~\ref{assumption1},~\ref{assumption2} and~\ref{assumption4}. We study the case where there exist reals $(a_{i,j})_{i,j\in\mathds N}$ such that for all $i,j\in\mathds N$, $h_{i,j}=a_{i,j} h$ for some $\pi$-canonical kernel $h:E^2 \to \mathds R$. For simplicity, we consider that $\mathds E_{\pi} h=0$ leading to $p_{i,j} = h_{i,j}$. Let us consider the specific example where $a_{i,j} = |j-i|^{-1}$ for $i\neq j$. In such setting, the coefficients $a_{i,j}$'s are weighting each summand in the U-statistic: the larger $|j-i|$, the smaller is the contribution of the term indexed by $(i,j)$ in the sum. As a result, interpreting indexes as time steps, the $a_{i,j}$'s can be understood as forgetting factors. Since
\begin{align*}
B_n^2&\leq A^2 \max\big\{\max_i \sum_{j=i+1}^n |j-i|^{-2} \; , \; \max_j \sum_{i=1}^{j-1} |j-i|^{-2}\big\} \leq A^2 \sum_{j=2}^n |j-1|^{-2}\leq A^2 \frac{\pi^2}{6},\\
\text{and} \quad C_n^2&\leq A^2 \sum_{j=2}^n \sum_{i=1}^{j-1} |j-i|^{-2}\leq A^2 \sum_{s=1}^n \frac{s}{s^2} \leq A^2 \big(1+\int_1^n\frac{1}{x}dx\big) \leq A^2(1+\log n),
\end{align*}
Theorem~\ref{mainthm3} ensures that there exist constants $\beta,\kappa>0$ such that for any $u\geq 1$ it holds with probability at least $1-\beta e^{-u}\log n$,
\[\frac{2}{n(n-1)}U_{\mathrm{stat}}(n)\leq \kappa A \log n \left(\log(n)\frac{\sqrt u}{n^{3/2}} + \left[\frac{u}{n}\right]^{3/2} +\left[\frac{u}{n}\right]^2\right). \]
In particular, with probability at least $1-\beta\frac{\log n }{n}$ we have $\frac{2}{n(n-1)} U_{\mathrm{stat}}(n) \leq 3\kappa A\frac{\log^{5/2} n}{n^{3/2}}$. This convergence rate improves significantly the one obtained from an Hoeffding-type concentration inequality like Eq.\eqref{eq:mainthm2} that would lead to $U_{\mathrm{stat}}(n) \leq 2\kappa A\frac{\log^{3/2} n}{n}$ with probability at least $1-\beta \frac{\log n}{n}$.

\section{Proofs}
\label{sec:proofs}

\subsection{Proofs of Theorems~\ref{mainthm} and~\ref{mainthm3}}

Our proof is inspired from~\cite{houdre2002} where a Bernstein-type inequality is shown for U-statistics of order~$2$ in the independent setting. Their proof relies on the \textit{canonical} property of the kernel functions which endowed the U-statistic with a martingale structure. We want to use a similar argument and we decompose~$ U_{\mathrm{stat}}(n)$ to recover the martingale property for each term (except for the last one). Considering for any~$l \geq 1$ the~$\sigma$-algebra~$G_l=\sigma(X_1, \dots,X_l)$, the notation~$\mathds E_{l}$ refers to the conditional expectation with respect to~$G_l$. Then we decompose~$ U_{\mathrm{stat}}(n)$ as follows,
\begin{align}  U_{\mathrm{stat}}(n)& =   M_{\mathrm{stat}}^{(t_n)}(n) +R_{\mathrm{stat}}^{(t_n)}(n),\label{decompo-ustat}\\
\text{with}\qquad  M_{\mathrm{stat}}^{(t_n)}(n)&=\sum_{k=1}^{t_n}  \sum_{i<j} \left(\mathds E_{j-k+1}[h_{i,j}(X_i,X_j)]-\mathds E_{j-k}[h_{i,j}(X_i,X_j)]\right),\notag \\
R_{\mathrm{stat}}^{(t_n)}(n)&=\sum_{i<j} \left(\mathds E_{j-t_n}[h_{i,j}(X_i,X_j)]-\mathds E \left[h_{i,j}(X_i,X_j)\right]\right),\notag 
\end{align}
and where~$t_n$ is an integer that scales logarithmically with~$n$. We recall that~$t_n:=\floor{r\log n}$ with~$r > 2\left(\log(1/\rho)\right)^{-1}$ where $\rho\in (0,1)$ is a constant characterizing the uniform ergodicity of the Markov chain (see Assumption~\ref{assumption1}). By convention, we assume here that for all~$k<1$,~$\mathds E_{k}[\cdot] := \mathds E[\cdot]$. Hence the first term that we will consider is given by
\[U_n = \sum_{1\leq i<j \leq n}  h^{(0)}_{i,j}(X_i,X_{j-1},X_j),\]
where for all ~$x,y,z \in E$, $h^{(0)}_{i,j}(x,y,z) = h_{i,j}(x,z)-  \int_w h_{i,j}(x,w)P(y,dw).$

We provide a detailed proof of a concentration result for~$U_n$ by taking advantage of its martingale structure. Reasoning by induction, we show that the~$t_n-1$ following terms involved in the decomposition~\eqref{decompo-ustat} of~$ U_{\mathrm{stat}}(n)$ can be handled using a similar approach. Since the last term $R_{\mathrm{stat}}^{(t_n)}(n)$ of the decomposition~\eqref{decompo-ustat} has not a martingale property, another argument is required. We deal with $R_{\mathrm{stat}}^{(t_n)}(n)$ exploiting the uniform ergodicity of the Markov chain~$(X_i)_{i \geq 1}$ which is guaranteed by Assumption~\ref{assumption1} (see~\cite[Theorem~8]{Roberts04}). 

The cornerstones of our approach are the following two propositions whose proofs are postponed to Section~\ref{proof:mainproposition} and Section~\ref{proof:mainproposition-remaining} respectively.

\begin{proposition}
\label{mainproposition} Let~$n \geq 2$. We keep the notations of Sections~\ref{sec:assumptions} and~\ref{subsec:csts}. We suppose Assumptions~\ref{assumption1},~\ref{assumption2} and~\ref{assumption4} described in Section~\ref{sec:assumptions}. 
There exist two constants~$\beta,\kappa>0$ such that for any~$u>0$, 

\begin{enumerate}[label=\alph*)]
\item if Assumption~\ref{assumption0}$.(i)$ is satisfied, it holds with probability at least~$1-\beta e^{-u}\log(n),$
\[M_{\mathrm{stat}}^{(t_n)}(n)\leq  \kappa \log(n)\Big( \left[A \sqrt n\log n\right]\sqrt u + \left[A+  B_n\sqrt n\right] u +\left[2A\sqrt n \right]u^{3/2} + A u^2\Big).\]
\item if Assumption~\ref{assumption0}$.(ii)$ is satisfied, it holds with probability at least~$1-\beta e^{-u}\log(n),$
\[M_{\mathrm{stat}}^{(t_n)}(n)\leq  \kappa \log(n)\Big( \left[\mathbf{C}_n+A \sqrt n\log n\right]\sqrt u + \left[A+  B_n\sqrt n\right] u +\left[2A\sqrt n \right]u^{3/2} + A u^2  \Big).\]
\end{enumerate}
\end{proposition}

\begin{proposition}
\label{mainproposition-remaining} Let~$n \geq 2$. We keep the notations of Sections~\ref{sec:assumptions} and~\ref{subsec:csts}. We suppose Assumptions~\ref{assumption1},~\ref{assumption2} and~\ref{assumption4}. Then
\begin{enumerate}[label=\alph*)]
\item $
R^{(t_n)}_{\mathrm{stat}}(n)\leq  A \left( 2L+ nt_n\right) .$
\item if the Markov chain $(X_i)_{i\geq1}$ is stationary, $R^{(t_n)}_{\mathrm{stat}}(n)\leq 2LA \left( 1+ t_n+ t_n^2\right) .$
\end{enumerate}
\end{proposition}

\subsubsection{Proof of Theorem~\ref{mainthm}}
\label{sec:proof-mainthm}

We suppose Assumptions~\ref{assumption1},~\ref{assumption2},~\ref{assumption4} and~\ref{assumption0}$.(i)$ (respectively~\ref{assumption0}$.(ii)$). From the decomposition~\eqref{decompo-ustat} coupled with Proposition~\ref{mainproposition}$.a)$ (respectively Proposition~\ref{mainproposition}$.b)$) and Proposition~\ref{mainproposition-remaining}$.a)$, the result of Theorem~\ref{mainthm}$.a)$ (respectively Theorem~\ref{mainthm}$.b)$) is straightforward.

\subsubsection{Proof of Theorem~\ref{mainthm3}}
\label{proof-mainthm3}

We suppose Assumptions~\ref{assumption1},~\ref{assumption2} and~\ref{assumption4}. We assume in addition that the Markov chain is stationary which implies in particular that Assumption~\ref{assumption0}$.(ii)$ holds. From the decomposition~\eqref{decompo-ustat} coupled with Proposition~\ref{mainproposition}$.b)$ and Proposition~\ref{mainproposition-remaining}$.b)$, the result of Theorem~\ref{mainthm3} is straightforward. Note that in case Assumption~\ref{assumption0}$.(i)$ holds, the quantity~$C_n$ (involved in the sub-Gaussian regime of the tail) can be removed from the inequality by simply using Proposition~\ref{mainproposition}$.a)$ rather than Proposition~\ref{mainproposition}$.b)$.

\subsection{Proof of Proposition~\ref{mainproposition}}

\label{proof:mainproposition}

Let us recall that Proposition~\ref{mainproposition} requires either a mild condition on the initial distribution of the Markov chain or the fact that the kernels~$h_{i,j}$ do not depend on~$i$ (see Assumption~\ref{assumption0}). One only needs to consider different Bernstein concentration inequalities for sums of functions of Markov chains to go from one result to the other. In this section, we give the proof of Proposition~\ref{mainproposition} in the case where Assumption~\ref{assumption0}$.(i)$ holds. We specify the part of the proof that should be changed to get the result when~$h_{i,j}$ may depend on both~$i$ and~$j$ and when Assumption~\ref{assumption0}$.(ii)$ holds. We make this easily identifiable using the symbol \raisebox{1.22ex}{\resizebox{!}{1ex}{ \dbend}}. 

\subsubsection{Concentration of the first term of the decomposition of the U-statistic}
\label{sec:firstterm}

\paragraph{Martingale structure of the U-statistic}

Defining~$Y_j = \sum_{i=1}^{j-1} h^{(0)}_{i,j}(X_i,X_{j-1},X_j)$,~$U_n$ can be written as~$U_n  = \sum_{j=2}^n Y_j.$
Since \[\mathds E_{j-1}[Y_j] =\mathds E[Y_j \; | \; X_1, \dots, X_{j-1}]=0, \]
we know that~$(U_k)_{k \geq 2}$ is a martingale relative to the~$\sigma$-algebras~$G_l$,~$l\geq2$. This martingale can be extended to~$n= 0$ and~$n= 1$ by taking~$U_0=U_1= 0$,~$G_0=\{\emptyset,E\}$,~$G_1=\sigma(X_1)$.
We will use the martingale structure of~$(U_n)_n$ through the following Lemma.

\begin{lemma} (cf.~\cite[Lemma 3.4.6]{ginenick})\\ \label{martingale-gine}Let~$(U_m , G_m ), m \in \mathds{N}$, be a martingale with respect to a filtration~$G_m$ such that~$U_0 = U_1 = 0$. For each~$m \geq  1$ and~$k \geq 2$, define the angle brackets~$A^k_m = A^k_m(U)$ of the martingale~$U$ by
\[A^k_m = \sum_{i=1}^m\mathds{E}_{i-1}[(U_i - U_{i-1} )^k  ]\]
(and note~$A^k_1 = 0$ for all~$k$). Suppose that for~$\alpha > 0$ and all~$i \geq 1$,~$\mathds{E}[e^{  \alpha|U_i -U_{i-1} |}] < \infty$. Then
\[ \left( \epsilon_m:= e^{ \alpha U_m - \sum_{k\geq2} \alpha^k A_m^k /k! }, G_m \right)\,,\ m \in \mathds{N},\]
is a supermartingale. In particular,~$\mathds{E}[\epsilon_m] \leq \mathds{E}[ \epsilon_1]=1$, so that, if~$A^k_m \leq w^k_m$ for constants ~$w^k_m \geq 0$ ; then
\[\mathds{E}[e^{\alpha U_m}] \leq e^{ \sum_{k\geq2} \alpha^k w_m^k /k!}.\]
\end{lemma}

We will also use the following convexity result several times.
\begin{lemma} \label{lemma:inequality} \cite[page 179]{ginenick}
For all~$\theta_1,\theta_2,\epsilon \geq 0$, and for all integer~$k\geq1$, \[(\theta_1+\theta_2)^k \leq (1+\epsilon)^{k-1}\theta_1^k+(1+\epsilon^{-1})^{k-1}\theta_2^k.\]
\end{lemma}

\noindent
For all~$k\geq 2$ and~$n\geq 1$, we have using Assumption~\ref{assumption4}:
{\small\begin{align*}
&A^k_n  = \sum_{j=2}^n \mathds{E}_{j-1} \Bigg[ \sum_{i=1}^{j-1} h^{(0)}_{i,j}(X_i,X_{j-1},X_j) \Bigg]^k\leq V^k_n :=\sum_{j=2}^n \mathds{E}_{j-1} \Bigg| \sum_{i=1}^{j-1} h^{(0)}_{i,j}(X_i,X_{j-1},X_j) \Bigg|^k \\
&=   \sum_{j=2}^n \mathds{E}_{j-1} \Bigg| \sum_{i=1}^{j-1}\big( h_{i,j}(X_i,X_j)-\mathds E_{\tilde X \sim \pi}[h_{i,j}(X_i,\tilde X)]+\mathds E_{\tilde X \sim \pi}[h_{i,j}(X_i,\tilde X)]-\mathds E_{j-1}[h_{i,j}(X_i,X_j)]  \big)\Bigg|^k\\
&=   \sum_{j=2}^n \mathds{E}_{j-1} \Bigg| \sum_{i=1}^{j-1} \left( p_{i,j}(X_i,X_j)+m_{i,j}(X_i,X_{j-1})\right) \Bigg|^k,
\end{align*}}
\[\text{where}\quad p_{i,j}(x,z) = h_{i,j}(x,z)- \mathds E_{ \pi}[h_{i,j}]\quad \text{ and } \quad m_{i,j}(x,y) =   \mathds E_{ \pi}[h_{i,j}] - \int_z h_{i,j}(x,z)P(y,dz).\]
Using Lemma~\ref{lemma:inequality} with~$\epsilon=1/2$, we deduce that
\begin{align*}
V^k_n& \leq \sum_{j=2}^n \mathds{E}_{j-1}\Bigg( \Bigg| \sum_{i=1}^{j-1} p_{i,j}(X_i,X_j) \Bigg|+  \Bigg| \sum_{i=1}^{j-1}m_{i,j}(X_i,X_{j-1}) \Bigg| \Bigg)^k\\
& \leq \left(\frac{3}{2}\right)^{k-1}\sum_{j=2}^n \mathds{E}_{j-1} \Bigg| \sum_{i=1}^{j-1} p_{i,j}(X_i,X_j)\Bigg|^k +3^{k-1}\sum_{j=2}^n \mathds{E}_{j-1} \Bigg| \sum_{i=1}^{j-1}m_{i,j}(X_i,X_{j-1}) \Bigg|^k.
\end{align*}

\noindent
Let us remark that

\begin{align*}
 \sum_{j=2}^n \mathds{E}_{j-1} \Bigg| \sum_{i=1}^{j-1}m_{i,j}(X_i,X_{j-1}) \Bigg|^k&= \sum_{j=2}^n \mathds{E}_{j-1} \Bigg| \sum_{i=1}^{j-1}\left(\mathds E_{\tilde X \sim \pi}[h_{i,j}(X_i,\tilde X)]-\mathds E_{j-1}[h_{i,j}(X_i,X_j)]\right) \Bigg|^k\\
&= \sum_{j=2}^n  \Bigg| \sum_{i=1}^{j-1}\left(\mathds E_{\tilde X \sim \pi}[h_{i,j}(X_i,\tilde X)]-\mathds E_{j-1}[h_{i,j}(X_i,X_j)]\right) \Bigg|^k\\
&= \sum_{j=2}^n  \Bigg|\mathds E_{j-1}\Bigg[ \sum_{i=1}^{j-1}\left(\mathds E_{\tilde X \sim \pi}[h_{i,j}(X_i,\tilde X)]-h_{i,j}(X_i,X_j)\right) \Bigg] \Bigg|^k\\
&\leq \sum_{j=2}^n  \mathds E_{j-1}\Bigg|\sum_{i=1}^{j-1}\left(\mathds E_{\tilde X \sim \pi}[h_{i,j}(X_i,\tilde X)]-h_{i,j}(X_i,X_j) \right) \Bigg|^k,
\end{align*}
where the last inequality comes from Jensen's inequality. We obtain the following upper-bound for~$V_n^k,$
\begin{align*}
V^k_n & \leq  2\times 3^{k-1}\sum_{j=2}^n \mathds{E}_{j-1} \Bigg| \sum_{i=1}^{j-1} p_{i,j}(X_i,X_j)\Bigg|^k\leq 2\times 3^{k-1}\delta_M\sum_{j=2}^n \mathds{E}_{X'_j} \Bigg| \sum_{i=1}^{j-1} p_{i,j}(X_i,X'_j)\Bigg|^k,
\end{align*}

\noindent
where the random variables~$(X'_j)_j$ are i.i.d. with distribution~$\nu$ (see Assumption~\ref{assumption2}).~$\mathds{E}_{X'_j}~$ denotes the expectation on the random variable~$X'_j$.

\begin{lemma} (cf.~\cite[Ex.1 Section 3.4]{ginenick})\label{duality}
Let~$Z_j$ be independent random variables with respective probability laws~$P_j$. Let~$k>1$, and consider functions~$f_1,\dots ,f_N$ where for all~$j\in [N]$,~$f_j\in L^k(P_j)$. Then the duality of~$L^p$ spaces and the independence of the variables~$Z_j$ imply that \[\left( \sum_{j=1}^N\mathds E\left[|f_j(Z_j)|^k\right]\right)^{1/k} = \underset{\sum_{j=1}^N\mathds E |\xi_j(Z_j)|^{k/(k-1)}=1}{\sup} \sum_{j=1}^N\mathds E\left[ f_j(Z_j)\xi_j(Z_j)\right],\]
where the sup runs over~$\xi_j \in L^{k/(k-1)}(P_j)$.
\end{lemma}

\noindent
Then by the duality result of Lemma~\ref{duality},
\begin{align*}
\left(V^k_n\right)^{1/k}&\leq \left( 2\delta_M\times 3^{k-1}  \sum_{j=2}^n \mathds{E}_{X'_j} \Bigg| \sum_{i=1}^{j-1} p_{i,j}(X_i,X'_j)\Bigg|^k\right)^{1/k} \\
& \leq \left(2\delta_M\right)^{1/k}\underset{\xi \in \mathcal{L}_k}{\sup}   \sum_{j=2}^n \sum_{i=1}^{j-1} \mathds{E}_{X'_j} \Bigg[  p_{i,j}(X_i,X'_j) \xi_j(X'_j) \Bigg] \\
\text{where} \;& \mathcal{L}_k = \Bigg\{ \xi = (\xi_2, \dots, \xi_{n}) \; s.t. \; \forall 2 \leq j \leq n, \; \xi_j \in L^{k/(k-1)}(\nu)   \; with \; \sum_{j=2}^n\mathds{E}
|\xi_j(X'_j)|^{\frac{k}{k-1}}=1 \Bigg\}.\\
&=  \left(2\delta_M\right)^{1/k}\underset{\xi \in \mathcal{L}_k}{\sup}    \sum_{i=1}^{n-1} \sum_{j=i+1}^{n}\mathds{E}_{X'_j} \Bigg[  p_{i,j}(X_i,X'_j) \xi_j(X'_j) \Bigg]
\end{align*}

Let us denote by~$F$ the subset of the set~$\mathcal F(E,\mathds R)$ of all measurable functions from~$(E,\Sigma)$ to~$(\mathds R,\mathcal B(\mathds R))$ that are bounded by~$A$. We set~$S:= E \times F^{n-1}$. For all~$i \in [n]$, we define~$W_i$ by
\[W_i:=\big(X_i,\underbrace{0,\dots,0}_{(i-1) \text{ times}},p_{i,i+1}(X_i,\cdot),p_{i,i+2}(X_i,\cdot), \dots, p_{i,n}(X_i,\cdot)\big)\in S.\]
Hence for all~$i \in [n]$,~$W_i$ is~$\sigma(X_i)$-measurable. We define for any~$\xi=(\xi_2, \dots, \xi_n) \in \prod_{i=2}^n L^{k/(k-1)}(\nu)$ the function \[\forall w=(x,p_2,\dots,p_n) \in S, \quad f_{\xi}(w) = \sum_{j=2}^n \int p_j(y) \xi_j(y) d\nu(y).\] Then setting~$\mathcal{F}=\{f_{\xi} : \sum_{j=2}^n \mathds E|\xi_j(X'_j)|^{k/(k-1)}=1\}$, we have
\[(V_n^k)^{1/k}\leq (2\delta_M)^{1/k}\underset{f_{\xi} \in \mathcal{F}}{\sup}  \sum_{i=1}^{n-1} f_{\xi}(W_i) .\]
By the separability of the~$L^p$ spaces of finite measures,~$\mathcal{F}$ can be replaced by a countable subset~$\mathcal{F}_0$.
To upper-bound the tail probabilities of~$U_n$, we will bound the variable~$V_n^k$ on sets of large probability using Talagrand's inequality. Then we will use Lemma~\ref{martingale-gine} on these sets by means of optional stopping.

\paragraph{Application of Talagrand's inequality for Markov chains}

The proof of Lemma~\ref{lemma:talagrand} is provided in Section \ref{proof-lemma-talagrand} of the Appendix (and relies mainly on the Talagrand concentration result from \cite{samson2000concentration}).

\begin{lemma}\label{lemma:talagrand}
Let us denote \[Z=\underset{f_{\xi} \in \mathcal{F}}{\sup} \sum_{i=1}^{n-1} f_{\xi}(W_i),\quad \sigma^2_k=\mathds E\left[ \sum_{i=1}^{n-1}\underset{f_{\xi} \in \mathcal{F}}{\sup} f_{\xi}(W_i)^2 \right] \quad \text{and} \quad b_k= \sup_{w \in  S} \sup_{f_{\xi} \in \mathcal F} |f_{\xi}(w)|.\] 
Then it holds for any~$t>0$,
\begin{align*}
&\mathds P \left( Z > \mathds E[Z]+t\right)\leq   \exp\left( -\frac{1}{8 \|\Gamma\|^2} \min \left( \frac{t^2}{4\sigma^2_k} \; , \; \frac{t}{b_k} \right) \right),
\end{align*}
where~$\Gamma$ is a~$n\times n$ matrix defined in Section \ref{proof-lemma-talagrand} %~\ref{proof-lemma-talagrand} 
which satisfies~$\|\Gamma\|\leq \frac{2L}{1-\rho}.$
\end{lemma}

Using Lemma~\ref{lemma:talagrand}, we deduce that for any~$t >0$,
\[\mathds P\left( (V_n^k)^{1/k} \geq (2\delta_M)^{1/k}  \mathds E[Z]+(2\delta_M)^{1/k} t\right) \leq  \exp \left( -\frac{1}{8 \|\Gamma\|^2} \min \left(  \frac{t^2}{4 \sigma^2_k}, \frac{t}{b_k}\right)\right),\]
which implies that for any~$x \geq 0$,
\[\mathds P\left( (V_n^k)^{1/k} \geq (2\delta_M)^{1/k}  \mathds E[Z]+(2\delta_M)^{1/k} 2\sigma_k \sqrt{x}+(2\delta_M)^{1/k} b_k x\right) \leq  \exp \left( -\frac{x}{8 \|\Gamma\|^2}\right).\]
Using the change of variable~$x=k8\|\Gamma\|^2u$ with~$u\geq0$ in the previous inequality leads to
\[\mathds P\left( \bigcup_{k=2}^{\infty}(V_n^k)^{1/k} \geq (2\delta_M)^{1/k} \mathds E[Z]+(2\delta_M)^{1/k} \sigma_k 3\|\Gamma\|\sqrt{ku}+(2\delta_M)^{1/k} k8\|\Gamma\|^2b_ku\right) \leq 1.62 e^{-u},\]
because \[1\wedge \sum_{k=2}^{\infty} \exp \left( -ku\right) \leq 1\wedge \frac{1}{e^u(e^u-1)}=\left(e^u \wedge \frac{1}{e^u-1}\right)e^{-u}\leq \frac{1+\sqrt{5}}{2}e^{-u}\leq 1.62e^{-u}.\]

Using Lemma~\ref{lemma:inequality} twice and using Holder inequality to bound~$b_k$ and~$\sigma_k^2$, we obtain~\eqref{vnk-wnk} from Lemma~\ref{lemma:bound-Vnk-proba}. The proof of Lemma~\ref{lemma:bound-Vnk-proba} is postponed to Section \ref{subsec:proof-bound-Vnk-proba}.
\begin{lemma}
\label{lemma:bound-Vnk-proba}
For any $u>0$, we denote \begin{align*}w^k_n&:= ( (1+\epsilon)^{k-1}2\delta_M \left(\mathds E[Z]\right)^k +2\delta_M(1+\epsilon^{-1})^{2k-2}\left(8\|\Gamma\|^2\right)^k(nA^2)A^{k-2}( ku)^k\\
& \quad +(1+\epsilon)^{k-1}(1+\epsilon^{-1})^{k-1}2\delta_M \left(3\|\Gamma\|\right)^k  \mathfrak B_0^2A^{k-2} (n ku)^{k/2},\end{align*}
\begin{equation} \label{eq:B0} \text{with} \; \; \mathfrak B_0^2:=  \max\big[ \max_{i} \big\|\sum_{j=i+1}^n \mathds E_{X \sim \nu}\left[ p_{i,j}^2(\cdot,X) \right] \big\|_{\infty} , \; \max_{j} \big\| \sum_{i=1}^{j-1} \mathds{E}_{X \sim \pi}[p_{i,j}^2(X,\cdot)] \big\|_{\infty}\big] \leq B_n^2,\end{equation}
where the dependence in~$u$ of~$w_n^k$ is leaved implicit. Then it holds
\begin{equation}\label{vnk-wnk} \mathds P \left( V_n^k \leq w_n^k \quad \forall k \geq2 \right) \geq 1-1.62e^{-u}.\end{equation}
\end{lemma}

\paragraph{Bounding $(\mathds E [Z])^k$.}\mbox{}\\
\label{proof:EZK}
\scalebox{0.95}{\dnote{The way we bound~$(\mathds E [Z])^k$ is the only part of the proof that needs to be modified to get the concentration result when Assumption~\ref{assumption0}$.(i)$ or Assumption~\ref{assumption0}$.(ii)$ holds. This is where we can use different Bernstein concentration inequalities according to whether the splitting method is applicable or not (see Section~\ref{sec:technical-assumption} for details). Here we present the approach when~$h_{i,j}\equiv h_{1,j}, \; \forall i,j$ (i.e. when Assumption~\ref{assumption0}.$(i)$ is satisfied). We refer to Section \ref{apdx:EZ-chi} of the Appendix for the details regarding the way we bound~$(\mathds E Z)^k$ when  Assumption~\ref{assumption0}$.(ii)$ holds.}}

Using Jensen inequality and Lemma~\ref{duality}, we obtain

\begin{align*}(\mathds E[Z])^k&\leq \mathds E[Z^k]=\mathds E\left[ \left( \underset{\xi \in \mathcal L_k}{\sup} \sum_{i=1}^{n-1}\sum_{j=i+1}^n\mathds E_{X'_j}[p_{i,j}(X_i,X'_j)\xi_j(X'_j)]  \right)^k\right]\\
&=\mathds E\left[  \sum_{j=2}^n\mathds E_{X'_j}\left[\left| \sum_{i=1}^{j-1}p_{i,j}(X_i,X'_j)\right|^k\right] \right]=  \sum_{j=2}^n\mathds E\left[\left| \sum_{i=1}^{j-1}p_{i,j}(X_i,X'_j)\right|^k\right] ,
\end{align*}
where we recall that~$\mathds{E}_{X'_j}~$ denotes the expectation on the random variable~$X'_j$. One can remark that conditionally to~$X'_j$, the quantity~$\sum_{i=1}^{j-1}p_{i,j}(X_i,X'_j)$ is a sum of function of the Markov chain~$(X_i)_{i\geq 1}$. Hence to control this term, we apply a Bernstein inequality for Markov chains.

Let us consider some~$j \in [n]$ and some~$x \in E$. Using the notations of Section~\ref{splitting}, we define
\[\forall l \in \{0, \dots,n\}, \quad Z_l^j(x)=\sum_{i=m(S_l+1)}^{m(S_{l+1}+1)-1}p_{i,j}(X_i,x).\]
By convention, we set~$p_{i,j}\equiv 0$ for any~$i\geq j$. Let us consider~$N_j = \sup \{i \in \mathbb N \;:\; mS_{i+1}+m-1\leq j-1\}$.
Then using twice Lemma~\ref{lemma:inequality}, we have
\begin{align}
&\big| \sum_{i=1}^{j-1}p_{i,j}(X_i,x)\big|^k=\big| \sum_{l=0}^{N_j}Z_l^j(x)+\sum_{i=m(S_{N_j}+1)}^{j-1}p_{i,j}(X_i,x)\big|^k \notag\\
&\leq \big(\frac32\big)^{k-1}\big| \sum_{l=1}^{N_j}Z_l^j(x)\big|^k+3^{k-1}\big|\sum_{i=m(S_{N_j}+1)}^{j-1}p_{i,j}(X_i,x)\big|^k \notag\\
&\leq \big(\frac94\big)^{k-1}\big|\sum_{l=0}^{\floor{N_j/2}}Z_{2l}^j(x)\big|^k+\big(\frac92\big)^{k-1}\big|\sum_{l=0}^{\floor{(N_j-1)/2}}Z_{2l+1}^j(x)\big|^k+3^{k-1}\big|\sum_{i=m(S_{N_j}+1)}^{j-1}p_{i,j}(X_i,x)\big|^k. \label{eq:bern-Tau}
\end{align}
\vspace{-0.45cm}

We have~$|\sum_{i=m(S_{N_j}+1)}^{j-1}p_{i,j}(X_i,x)|\leq AmT_{N_j+1}$.  So using the definition of the Orlicz norm and the fact that the random variables~$(T_i)_{i\geq 2}$ are i.i.d., it holds for any~$t\geq 0$,
\begin{align*}
\mathds P\left(\left|\sum_{i=m(S_{N_j}+1)}^{j-1}p_{i,j}(X_i,x)\right|\geq t\right)&\leq \mathds P(T_{N_j+1}\geq \frac{t}{Am})\leq \mathds P(\max(T_1,T_2)\geq \frac{t}{Am})\\
&\leq \mathds P(T_1\geq \frac{t}{Am})+\mathds P(T_2\geq \frac{t}{Am})\leq 4 \exp(-\frac{t}{Am\tau}).
\end{align*}
Hence, using that for an exponential random variable $G$ with parameter~$1$, $\mathds E[G^p] = p! \;\; \forall p \geq 0$,
\begin{align*}
&\mathds E \left[\left|\sum_{i=m(S_{N_j}+1)}^{j-1}p_{i,j}(X_i,x)\right|^k\right] =4 \int_0^{+\infty } \mathds P\left(\left|\sum_{i=m(S_{N_j}+1)}^{j-1}p_{i,j}(X_i,x)\right|^k\geq t\right)dt\\
&\leq  4 \int_0^{+\infty}\exp(-\frac{t^{1/k}}{Am\tau})\leq  4(Am\tau)^k \int_0^{+\infty}\exp(-v)kv^{k-1}dv=4(Am\tau)^k k!,
\end{align*}

The random variable~$Z_{2l}^j(x)$ is~$\sigma(X_{m(S_{2l}+1)},\dots,X_{m(S_{2l+1}+1)-1})$-measurable. Let us insist that this holds because we consider that~$h_{i,j}\equiv h_{1,j}, \; \forall i,j$ which implies that~$p_{i,j} \equiv p_{1,j},\; \forall i,j$. Hence for any~$x\in E$, the random variables~$(Z_{2l}^j(x))_l$ are independent (see Section~\ref{splitting}). Moreover, one has that for any~$l$,~$\mathds E[Z^j_{2l}(x) ]=0$. This is due to~\cite[Eq.(17.23) Theorem~17.3.1]{tweedie} together with Assumption~\ref{assumption4} which gives that 
~$ \forall x' \in E, \quad   \mathds E_{X \sim \pi} [p_{i,j}(X,x') ]=0.$
Let us finally notice that for any~$x\in E$ and any~$l\geq0$,~$|Z_{2l}^j(x)| \leq Am  T_{2l+1}$, so~$\|Z_{2l}^j(x)\|_{\psi_1}\leq Am\max(\|T_1\|_{\psi_1},\|T_2\|_{\psi_1})\leq Am \tau$. 
First, we use Lemma~\ref{lemma:decoupling} to obtain that 
\[\mathds E\left|\sum_{l=0}^{\floor{N_j/2}}Z_{2l}^j(x)\right|^k \leq \mathds E\max_{0\leq s \leq n-1}\left|\sum_{l=0}^{s}Z_{2l}^j(x)\right|^k\leq 2\times 4^{k} \mathds E \left|\sum_{l=0}^{n-1}Z_{2l}^j(x)\right|^k,\]
where for the last inequality we gathered~\eqref{decoupling:2} with the left hand side of~\eqref{decoupling:1} from Lemma~\ref{lemma:decoupling}.

\begin{lemma} (cf.~\cite[Lemma 1.2.6]{gine-decoupling})\\ \label{lemma:decoupling}
Let us consider some separable Banach space~$B$ endowed with the norm~$\|\cdot\|$. Let~$X_i$,~$i\leq n$, be independent centered~$B$-valued random variables with norms~$L_p$ for some~$p\geq 1$ and let~$\epsilon_i$ be independent Rademacher random variables independent of the variables~$X_i$. Then
\begin{align}\label{decoupling:1}2^{-p}\mathds E \left\| \sum_{i=1}^n \epsilon_i X_i\right\|^p & \leq \mathds E \left\| \sum_{i=1}^n  X_i\right\|^p\leq 2^{p}\mathds E \left\| \sum_{i=1}^n \epsilon_i X_i\right\|^p,\\
\text{and} \quad\label{decoupling:2}\mathds E \max_{k \leq n}\left\| \sum_{i=1}^k  X_i\right\|^p &\leq 2^{p+1}\mathds E \left\| \sum_{i=1}^n \epsilon_i X_i\right\|^p\end{align}
\end{lemma}
Similarly, the random variables~$(Z_{2l+1}^j(x))_l$ are independent and satisfy for any~$l$,~$\mathds E[Z^j_{2l+1}(x) ]=0$. With an analogous approach, we get that
\[\mathds E\left|\sum_{l=0}^{\floor{(N_j-1)/2}}Z_{2l+1}^j(x)\right|^k \leq \mathds E\max_{0\leq s \leq n-1}\left|\sum_{l=0}^{s}Z_{2l+1}^j(x)\right|^k\leq 2\times 4^{k} \mathds E \left|\sum_{l=0}^{n-1}Z_{2l+1}^j(x)\right|^k.\]
Let us denote for any~$j \in [n]$,~$\mathds E_{|X'_j}$ the conditional expectation with respect to the~$\sigma$-algebra~$\sigma(X'_j)$. Coming back to~\eqref{eq:bern-Tau}, we proved that
\begingroup
 \allowdisplaybreaks
\begin{align}
&\mathds E_{|X'_j}\left| \sum_{i=1}^{j-1}p_{i,j}(X_i,X'_j)\right|^k \leq \left(\frac94\right)^{k-1}\mathds E_{|X'_j}\left|\sum_{l=0}^{\floor{N_j/2}}Z_{2l}^j(X'_j)\right|^k\notag\\
&+\left(\frac92\right)^{k-1}\mathds E_{|X'_j}\left|\sum_{l=0}^{\floor{(N_j-1)/2}}Z_{2l+1}^j(X'_j)\right|^k+3^{k-1}\mathds E_{|X'_j}\left|\sum_{i=m(S_{N_j}+1)}^{j-1}p_{i,j}(X_i,X'_j)\right|^k \notag\\
&\leq 2 \times  9^{k}\mathds E_{|X'_j} \left|\sum_{l=0}^{n-1}Z_{2l+1}^j(X'_j)\right|^k   +2\times 18^{k}   \mathds E_{|X'_j} \left|\sum_{l=0}^{n-1}Z_{2l}^j(X'_j)\right|^k + 4 (3Am\tau)^k k!.\label{eq:bern-tau-2}
\end{align}
\endgroup
It remains to bound the two expectations in~\eqref{eq:bern-tau-2}. The two latter expectations will be controlled similarly and we give the details for the first one. We use the following Bernstein's inequality with the sequence of random variables~$(Z_{2l+1}^j(x))_l$.
\begin{lemma} (Bernstein's~$\psi_1$ inequality,~\cite[Lemma 2.2.11]{VW00} and the subsequent remark).\label{bernstein-psi1} \\If~$Y_1,\dots,Y_n$ are independent random variables such that~$\mathds EY_i= 0$ and~$\|Y_i\|_{\psi_1}\leq \tau$, then for every~$t >0,$
\[ \mathds P \big( \big|\sum_{i=1}^n Y_i\big|>t\big)\leq 2\exp \big( -\frac{1}{K} \min \big( \frac{t^2}{n\tau^2},\frac{t}{\tau}\big) \big) ,\]
for some universal constant~$K>0$ ($K=8$ fits).
\end{lemma}
We obtain 
\[ \mathds P \left(\left|\sum_{l=0}^{n-1}Z_{2l+1}^j(x)\right|>t\right)\leq 2\exp \left( -\frac{1}{K} \min \left( \frac{t^2}{nA^2m^2\tau^2},\frac{t}{Am\tau}\right) \right) .\]
We deduce that for any~$x \in E$, any~$j \in [n]$ and any~$t\geq 0$,
\begingroup
 \allowdisplaybreaks
\begin{align*}
\mathds E\left[  \left|\sum_{l=0}^{n-1} Z^j_{2l+1}(x)\right|^k\right]& = \int_0^{\infty}  \mathds P\left(\left|\sum_{l=0}^{n-1} Z^j_{2l+1}(x)\right|^k>t\right)dt\\
&= 2 \int_0^{\infty} \exp \left(- \frac{1}{K} \min \left( \frac{t^{2/k}}{nA^2m^2 \tau^2},\frac{t^{1/k}}{Am\tau}\right)  \right) dt.
\end{align*}
\endgroup
Let us remark that $ \frac{t^{2/k}}{A^2m^2 n\tau^2} \leq \frac{t^{1/k}}{Am\tau} \Leftrightarrow t \leq (nA\tau m)^k.$
Hence for any~$j \in [n]$,
\begingroup
\allowdisplaybreaks
\begin{align*}
&\mathds E\left[  \left|\sum_{l=0}^{n-1} Z^j_{2l+1}(X'_j)\right|^k\right]\\
&\leq  2 \int_0^{(nA\tau m)^k} \exp \left(-   \frac{t^{2/k}}{KnA^2m^2\tau^2} \right) dt +2 \int_0^{\infty} \exp\left( - \frac{t^{1/k}}{KAm\tau}\right)dt.\\
&\leq 2 \int_0^{n/K} \exp \left(-   v \right) \frac{k}{2} v^{k/2-1}  \left( \sqrt{K}n^{1/2}  A\tau m \right)^k dv +2 \int_0^{\infty} \exp\left( -v\right) k v^{k-1} (KAm \tau)^kdv.\\
&\leq 2\int_0^{n/K} \exp \left(-   v \right) \frac{k}{2} v^{k/2-1}  \left( \sqrt{K}n^{1/2}  A\tau m \right)^k dv + 2 k \times (k-1)!   (KAm\tau)^k\\
&\leq k \left( \sqrt{K}n^{1/2}  A\tau m \right)^k\int_0^{n/K} \exp \left(-   v \right)  v^{k/2-1}   dv +  2k! (KAm\tau)^k ,
\end{align*}
\endgroup
where we used again that if~$G$ is an exponential random variable with parameter~$1$, then for any~$p \in \mathds N$,~$\mathds E[G^p] = p!$. Since for any real~$l \geq1$,
\begin{align*}
\int_0^{\frac{n}{K}} e^{-v} v^{l-1} dv&= \sum_{r=0}^{+\infty} \frac{(-1)^r}{r!}\int_0^{\frac{n}{K}} v^{r+l-1} dv= \sum_{r=0}^{+\infty} \frac{(-1)^r(\frac{n}{K})^{r+l}}{r!(r+l)} \leq \sum_{r=0}^{+\infty} \frac{(-1)^r(\frac{n}{K})^{r+l}}{r!\times l} \leq \frac{(\frac{n}{K})^{l}}{l}e^{-\frac{n}{K}},
\end{align*}
we get that
\[k \left( \sqrt{K}n^{1/2}  A\tau m \right)^k\int_0^{\frac{n}{K}}e^{-v}  v^{k/2-1}  dv \leq 2 \left( \sqrt{K}n^{1/2}A\tau m \right)^k e^{-\frac{n}{K}}\left(\frac{n}{K}\right)^{k/2} = 2 \left( nA\tau m \right)^k e^{-\frac{n}{K}}.\]
Hence we proved that for some universal constant~$K>1$,
\begin{align*}
&\mathds E\left[  \left|\sum_{l=0}^{n-1} Z^j_{2l+1}(x)\right|^k\right]\leq  2 \left( nA\tau m \right)^k e^{-n/K}+ 2k! (KAm \tau)^k\leq 4k! (KAm\tau)^k,
\end{align*}
since for all~$k \geq 2$,~$e^{-n/K}(n/K)^k /(k!)  \leq 1.$
Using a similar approach, one can show the same bound for the second expectation in~\eqref{eq:bern-tau-2}. We proved that for some universal constant~$K>1$,
\begin{align}
&(\mathds E[Z])^k \leq \sum_{j=2}^n\mathds E \left[ \mathds E_{|X'_j}\left| \sum_{i=1}^{j-1}p_{i,j}(X_i,X'_j)\right|^k\right] \leq 2 \times  9^{k}\sum_{j=2}^n\mathds E\left[\mathds E_{|X'_j} \left|\sum_{l=0}^{n-1}Z_{2l+1}^j(X'_j)\right|^k \right]  \notag\\
&\qquad \qquad +2\times 18^{k} \sum_{j=2}^n\mathds E\left[  \mathds E_{|X'_j} \left|\sum_{l=0}^{n-1}Z_{2l}^j(X'_j)\right|^k \right]+ 4\sum_{j=2}^n (3Am\tau)^k k! \notag\\
&\leq 2 n\times  18^{k} \times 4k! (KAm\tau)^k + 4n (3Am\tau)^k k!= 16n\times k! (KAm\tau)^k,\label{eq:boundEZk}
\end{align}
where in the last inequality, we still call~$K$ the universal constant defined by~$18K.$

\paragraph{Upper-bounding $U_n$ using the martingale structure}

Let \[T+1:=\inf \{l \in \mathds N : V_l^k \geq w_n^k \text{ for some }k\geq2\}.\]
Then, the event~$\{T\leq l\}$ depends only on~$X_1,\dots,X_l$ for all~$l\geq1$. Hence,~$T$ is a stopping time for the filtration~$(\mathcal{G}_l)_l~$ where~$\mathcal{G}_l=\sigma((X_i)_{i \in [l]})$ and we deduce that~$U_l^T:=U_{l\wedge T}$ for~$l=0, \dots ,n$ is a martingale with respect to~$(\mathcal{G}_l)_l$ with~$U^T_0=U_0= 0$ and~$U^T_1=U_1=0$. We remark that~$U^T_j-U^T_{j-1}=U_j-U_{j-1}$ if~$T\geq j$ and zero otherwise, and that~$\{T \geq j\}$ is~$\mathcal{G}_{j-1}$ measurable.
Then, the angle brackets of this martingale admit the following bound:
\begingroup
 \allowdisplaybreaks
\begin{align*}
A_n^k(U^T)&=\sum_{j=2}^{n} \mathds E_{j-1}[(U^T_j-U^T_{j-1})^k]\\
&\leq\sum_{j=2}^{n} \mathds E_{j-1}|U_j-U_{j-1}|^k\mathds 1_{T \geq j}= \sum_{j=2}^{n} \mathds E_{j-1}\left|\sum_{i=1}^{j-1}h^{(0)}(X_i,X_{j-1},X_j)\right|^k\mathds 1_{T \geq j}\\
&= \sum_{j=2}^{n-1} V_j^k \mathds 1_{T =j}+V_n^k \mathds 1_{T \geq n}\leq w_n^k\left( \sum_{j=2}^{n-1} \mathds 1_{T =j}+\mathds 1_{T \geq n}\right)\leq w_n^k,
\end{align*}
\endgroup
since, by definition of~$T$,~$V^k_j\leq w^k_n$ for all~$k$ on~$\{T\geq j\}$. Hence, Lemma~\ref{martingale-gine} applied to the martingale~$U^T_n$ implies
\[\mathds E e^{\alpha U^T_n} \leq \exp \left( \sum_{k\geq2}\frac{\alpha^k}{k!}w_n^k\right).\]
Also, since~$V^k_n$ is nondecreasing in~$n$ for each~$k$, inequality~\eqref{vnk-wnk} implies that
\begin{equation*} \mathds P(T<n) \leq \mathds P \left( V_n^k \geq w_n^k \quad \text{ for some } k \geq2 \right) \leq 1.62 e^{-u}.\end{equation*}
\noindent
Thus we deduce that for all~$s\geq0$,
\begin{equation}\label{ustats:final-bound}\mathds P(U_n\geq s) \leq \mathds P(U^T_n \geq s,T\geq n ) +\mathds P(T<n) \leq e^{-\alpha s} \exp\left( \sum_{k \geq2} \frac{\alpha^k}{k!} w_n^k \right) + 1.62 e^{-u}.\end{equation}
The final step of the proof consists in simplifying~$\exp \left( \sum_{k\geq2}\frac{\alpha^k}{k!}w_n^k\right)$.

\begin{align}
\sum_{k\geq2}\frac{\alpha^k}{k!}w_n^k &=  2\delta_M \sum_{k\geq2}\frac{\alpha^k}{k!}(1+\epsilon)^{k-1} (\mathds E[Z])^k \notag\\
&\qquad + 2\delta_M \sum_{k\geq2}\frac{\alpha^k}{k!}(2+\epsilon+\epsilon^{-1})^{k-1}\left(3\|\Gamma\|\right)^k    \mathfrak B_0^2A^{k-2} (n ku)^{k/2}\notag\\
&\qquad +2\delta_M\sum_{k\geq2}\frac{\alpha^k}{k!}(1+\epsilon^{-1})^{2k-2}\left(8\|\Gamma\|^2\right)^k  (nA^2)A^{k-2}( ku)^k=:a_1+a_2+a_3.\label{eq:defai}
\end{align}

\noindent Using the bound Eq.\eqref{eq:boundEZk} obtained on $\left(\mathds EZ\right)^k$, Lemma~\ref{lemma:bound-a1a2a3} bounds the three sums~$a_1, a_2$ and~$a_3.$
\begin{lemma}\label{lemma:bound-a1a2a3}
$\exp \left( \sum_{k\geq2}\frac{\alpha^k}{k!}w_n^k\right) \leq \exp\left( \frac{\alpha^2 W^2}{1-\alpha c}\right)$ where
\begin{align*} W &= 6\sqrt{\delta_M}(1+\epsilon)^{1/2}n^{1/2} KA\tau m\\
&\qquad +\sqrt{2\delta_M}(2+\epsilon+\epsilon^{-1})^{1/2}3\|\Gamma\|  \mathfrak B_0 \sqrt{nu}+\sqrt{2 \delta_M}A(1+\epsilon^{-1})8\|\Gamma\|^2 \sqrt{n}eu, \\
\text{and} \quad c&=\max\Bigg[ (1+\epsilon) KA\tau m\; , \;    (2+\epsilon+\epsilon^{-1})\left(3\|\Gamma\|\right)  A (n u)^{1/2} \;, \; (1+\epsilon^{-1})^{2}\left(8\|\Gamma\|^2\right)  A eu\Bigg].
\end{align*}
\end{lemma}

\noindent
The proof of Lemma~\ref{lemma:bound-a1a2a3} can be found in Section \ref{apdx:bound-a1a2a3} Using the result from Lemma~\ref{lemma:bound-a1a2a3} in~\eqref{ustats:final-bound} and taking~$s=2W\sqrt{u}+cu$ and~$\alpha = \sqrt{u}/(W+c\sqrt{u})$ in this inequality yields \[ \mathds P \left( U_n \geq 2 W \sqrt{u} + cu\right) \leq e^{-u}+1.62e^{-u}\leq (1+e)e^{-u} .\]

\noindent
By taking~$\epsilon=1/2$, we deduce that for any~$u\geq0$, it holds with probability at least~$1-(1+e)e^{-u}$
\begin{align*}  \sum_{i<j}h^{(0)}_{j}&(X_i,X_{j-1},X_j)
\leq \quad  12\sqrt{\delta_M}KA\tau m  \sqrt{nu}+18\sqrt{\delta_M}\|\Gamma\|  \mathfrak B_0\sqrt{ n}u \\
&+100\sqrt{\delta_M}\|\Gamma\|^2A\sqrt{n}eu^{3/2} +3KA\tau m  u +27A\|\Gamma\| \sqrt{n}u^{3/2} + 72A\|\Gamma\|^2 eu^2, \end{align*}

\noindent
Denoting~$ \kappa := \max\left( 12\sqrt{\delta_M}K\tau m  \; , \; 18\sqrt{\delta_M}\|\Gamma\|,   100\sqrt{\delta_M}\|\Gamma\|^2e,3K\tau m,72\|\Gamma\|^2 e\right),$
we have with probability at least~$1-(1+e)e^{-u}$
\[\sum_{i<j}h^{(0)}_{j}(X_i,X_{j-1},X_j)\leq \kappa \left( A \sqrt{n}  \sqrt{u} + (A  +  \mathfrak B_0\sqrt n) u +2A\sqrt n u^{3/2} + A  u^2 \right).\]

\subsubsection{Reasoning by descending induction with a logarithmic depth}
\label{desc-induction}

As previously explained, we apply a proof similar to the one of the previous subsection on the~$t_n:=\floor{r\log n}$ terms in the sum $M_{\mathrm{stat}}^{(t_n)}(n)$ (see~\eqref{decompo-ustat}), with~$r >2 \left(\log(1/\rho)\right)^{-1}$. Let us give the key elements to justify such approach by considering the second term of the sum $M_{\mathrm{stat}}^{(t_n)}(n)$, namely \begin{align}&\sum_{i<j} \left( \mathds E_{j-1}\left[ h_{i,j}(X_i,X_j)\right] - \mathds E_{j-2}\left[ h_{i,j}(X_i,X_j)\right] \right)\label{2nd-term-ustat}\\
= \quad & \underbrace{\sum_{i=1}^{n-2} \sum_{j=i+2}^{n} h^{(1)}_{i,j}(X_i,X_{j-2},X_{j-1})}_{=:U_{n-1}^{(1)}} + \underbrace{\sum_{i=1}^{n-1}\left\{\mathds E_{i}\left[ h_{i,i+1}(X_i,X_{i+1})\right] - \mathds E_{i-1}\left[ h_{i,i+1}(X_i,X_{i+1})\right]\right\}}_{=:(*)}\notag\end{align}
where $h^{(1)}_{i,j}(x,y,z) = \int_w h_{i,j}(x,w)P(z,dw) -\int_w h_{i,j}(x,w)P^2(y,dw).$ Using McDiarmid’s inequality for Markov chain (see~\cite[Corollary 2.10 and Remark 2.11]{paulin15}), we obtain Lemma~\ref{lemma:MCDIARMID}.

\begin{lemma} \label{lemma:MCDIARMID}
Let us consider $l\in \{1,\dots,t_n\}$.
 For any~$u>0$, it holds with probability at least~$1-2e^{-u}$,
\[\left| \sum_{i=1}^{n-1} \sum_{j=i+1}^{(i+l)\wedge n} \left( \mathds E_{j-l}\left[ h_{i,j}(X_i,X_j)\right] - \mathds E_{j-l-1}\left[ h_{i,j}(X_i,X_j)\right] \right) \right|\leq 3A t_n \sqrt{t_{mix}nu},\]
where~$t_{mix}$ is the mixing time of the Markov chain and is given by
\[t_{mix}:= \min \left\{t \geq 0 \; : \; \sup_{x} \|P^t(x, \cdot)- \pi\|_{\mathrm{TV}}<\frac14\right\}.\]

\end{lemma}

Lemma~\ref{lemma:MCDIARMID} allows to bound $(*)$ in~\eqref{2nd-term-ustat} (by choosing $l=1$). Now we aim at proving a concentration result for the term
\[U_{n-1}^{(1)}=\sum_{j=2}^{n-1} \sum_{i=1}^{j-1} h^{(1)}_{i,j}(X_i,X_{j-1},X_{j}),\]
using an approach similar to the one of the previous subsection. We state here the concentration result for $U^{(1)}_{n-1}$ with Lemma~\ref{lemma:induction} and a proof can be found in Section \ref{apdx:proof-lemma-induction}.

\begin{lemma}
\label{lemma:induction}For any~$u>0$, it holds with probability at least~$1-(1+e)e^{-u}$,
\begin{equation*}
 U^{(1)}_{n-1}=\sum_{i=1}^{n-2} \sum_{j=i+1}^{n-1} h^{(1)}_{i,j}(X_i,X_{j-1},X_{j})  \leq  \kappa \left( A \sqrt{n}  \sqrt{u} + (A  +  B_n\sqrt n) u +2A\sqrt n u^{3/2} + A  u^2\right)
\end{equation*}
\end{lemma}
Going back to~\eqref{2nd-term-ustat} and using both Lemmas~\ref{lemma:MCDIARMID} and~\ref{lemma:induction}, we get that for any~$u>0$, it holds with probability at least~$1-(1+e+2)e^{-u}$,
 \begin{align}& \sum_{i<j} \left( \mathds E_{j-1}\left[ h_{i,j}(X_i,X_j)\right] - \mathds E_{j-2}\left[ h_{i,j}(X_i,X_j)\right] \right)\nonumber\\
 \leq \quad & \kappa \left( A \sqrt{n}  \sqrt{u} + (A  +  B_n\sqrt n) u +2A\sqrt n u^{3/2} + A  u^2 \right)+3A t_n \sqrt{t_{mix}nu} \label{end-induction-1st}
 \end{align}

One can do the same analysis for the~$t_n$ first terms in the decomposition~\eqref{decompo-ustat}. Still denoting $\kappa$ the constant $\kappa+3\sqrt{t_{mix}}$, we get that for any~$u>0$ it holds with probability at least~$1-(3+e)e^{-u}t_n,$

\begin{equation*}M^{(t_n)}_{\mathrm{stat}}(n)\leq \kappa  t_n\left(  A t_n \sqrt{n}  \sqrt{u} + (A  +  B_n\sqrt n) u +2A\sqrt n u^{3/2} + A  u^2\right)
.\end{equation*}

\subsection{Proof of Proposition~\ref{mainproposition-remaining} }

\label{proof:mainproposition-remaining}
In the following, we assume that~$t_n\leq n$, otherwise $R^{(t_n)}_{\mathrm{stat}}(n)$ is an empty sum. Using our convention which states that for all~$k<1$,~$\mathds E_{k}[\cdot] := \mathds E[\cdot]$,
we need to control
\begin{align}
\left| R^{(t_n)}_{\mathrm{stat}}(n)\right|&=\left|\sum_{i<j} \left( \mathds E_{j-t_n}\left[ h_{i,j}(X_i,X_j)\right]-\mathds E\left[h_{i,j}( X_i,  X_j)\right]  \right)\right|\leq  (1)+(2),\label{eq:Rstat}
\end{align}
with denoting~$H_{i,j} = \mathds E_{j-t_n}\left[ h_{i,j}(X_i,X_j)\right]-\mathds E\left[h_{i,j}( X_i,  X_j)\right]~$,
\begin{align*}
(1)&:= \left| \sum_{i=1}^{n-t_n} \sum_{j=i+t_n}^n  H_{i,j} \right|= \left| \sum_{j=t_n+1}^{n} \sum_{i=1}^{j-t_n}  H_{i,j} \right| \\\:\text{and}\quad  (2)&:= \left| \sum_{i=1}^{n-1} \sum_{j=i+1}^{(i+t_n-1)\wedge n}  H_{i,j}\right|= \left| \sum_{j=2}^{n} \; \sum_{i=(j-t_n+1)\vee 1}^{j-1} H_{i,j} \right|.
\end{align*}
We start by bounding the term $(1)$ regardless of the initial distribution of the chain. We will bound in different ways the term $(2)$ depending on whether the Markov chain is stationary or not. Let us first bound the term~$(1)$ splitting it into two terms,
\begin{align*}
(1)=&\left| \sum_{j=t_n+1}^{n} \sum_{i=1}^{j-t_n}  \mathds E_{j-t_n}\left[ h_{i,j}(X_i,X_j)\right]-\mathds E\left[h_{i,j}( X_i,  X_j)\right] \right|\leq (1a) + (1b).
\end{align*}
Using Assumption~\ref{assumption4}, it holds~$\mathds E_{\pi}[h_{i,j}] = \mathds E_{\tilde X \sim \pi}[h_{i,j}(X_i,\tilde X)]=\int_{x} h_{i,j}(X_i,x) d\pi(x)$. Hence we get that
\begin{align*}
(1a):=&\left| \sum_{j=t_n+1}^{n} \sum_{i=1}^{j-t_n}  \mathds E_{j-t_n}\left[ h_{i,j}(X_i,X_j)\right]-\mathds E_{\pi}[h_{i,j}] \right|\\
&\leq \sum_{j=t_n+1}^{n} \left| \int_{x_j} \sum_{i=1}^{j-t_n}h_{i,j}(X_i,x_j) \left( P^{t_n}(X_{j-t_n},dx_j)- d\pi(x_j)\right) \right| \\
&\leq \sum_{j=t_n+1}^{n} \sup_{x_j} \left| \sum_{i=1}^{j-t_n}h_{i,j}(X_i,x_j)\right| \sup_{z}\| P^{t_n}(z,\cdot)- \pi\|_{\mathrm{TV}}\\
&\leq \sum_{j=t_n+1}^{n} \sup_{x_j} \left| \sum_{i=1}^{j-t_n}h_{i,j}(X_i,x_j)\right| L \rho^{t_n}\leq  \sum_{j=t_n+1}^{n} \sup_{x_j} \left| \sum_{i=1}^{j-t_n}h_{i,j}(X_i,x_j)\right| L \frac{1}{n^2}\leq LA,
\end{align*}
where in the penultimate inequality we used that~$\rho^{t_n}\leq\rho^{r\log(n)} = n^{r \log(\rho)}\leq n^{-2}~$. Indeed ~$2+r\log(\rho)<0$ because we choose~$r~$ such that~$r>2(\log(1/\rho))^{-1}.$ 

Using again Assumption~\ref{assumption4}, it holds~$\mathds E_{\pi}[h_{i,j}] =\int_{x_i} \chi P^i(dx_i) \int_x h_{i,j}(x_i,x)d\pi(x)$ where~$\chi$ is the initial distribution of the Markov chain~$(X_i)_{i\geq 1}$. We get that
\begin{align*}
(1b):=&\left| \sum_{j=t_n+1}^{n} \sum_{i=1}^{j-t_n} \mathds E_{\pi}[h_{i,j}]-\mathds E\left[h_{i,j}( X_i,  X_j)\right] \right|\\
&\leq \sum_{j=t_n+1}^{n} \sum_{i=1}^{j-t_n}\left| \int_{x_i} \int_{x_j} h_{i,j}(x_i,x_j) \chi P^i( dx_i)\left( P^{j-i}(x_i,dx_j)- d\pi(x_j)\right) \right| \\
&\leq \sum_{j=t_n+1}^{n} \sum_{i=1}^{j-t_n}  \|h_{i,j}\|_{\infty} \underbrace{\int_{x_i} \chi P^i( dx_i)}_{=1} \sup_{z} \int_{x_j}\left| P^{j-i}(z,dx_j)- d\pi(x_j)\right| \\
&\leq \sum_{j=t_n+1}^{n}  \sum_{i=1}^{j-t_n}\|h_{i,j}\|_{\infty}  L \rho^{j-i}\leq \sum_{j=t_n+1}^{n}  \sum_{i=1}^{j-t_n}\|h_{i,j}\|_{\infty}  L \rho^{t_n}\leq LA,
\end{align*}
where in the penultimate inequality we used that~$\rho^{t_n}\leq\rho^{r\log(n)} = n^{r \log(\rho)}\leq n^{-2}~$.

\subsubsection{Bounding (2) without stationarity}
Without assuming that the Markov chain is stationary, we bound  coarsely~$(2)$ as follows
\begin{align*}
(2)&=\left| \sum_{j=2}^{n} \; \sum_{i=(j-t_n+1)\vee 1}^{j-1}  \mathds E_{j-t_n}\left[ h_{i,j}(X_i,X_j)\right]-\mathds E\left[h_{i,j}( X_i, X_j)\right] \right|\leq Ant_n.
\end{align*}

This concludes the proof of Proposition~\ref{mainproposition-remaining}$.a)$ since we obtain, $R^{(t_n)}_{\mathrm{stat}}(n)\leq  A \left( 2L+ nt_n\right).$

\subsubsection{Bounding (2) with stationarity}

Considering now that the chain is stationary, we split~$(2)$ into three different contributions.
\begin{align*}
(2)&=\left| \sum_{j=2}^{n} \; \sum_{i=(j-t_n+1)\vee 1}^{j-1}  \mathds E_{j-t_n}\left[ h_{i,j}(X_i,X_j)\right]-\mathds E\left[h_{i,j}( X_i, X_j)\right] \right|\leq (2a) + (2b)+ (2c),\\
\text{with} \quad (2a)&:= \left| \sum_{j=2}^{n} \; \sum_{i=(j-t_n+1)\vee 1}^{j-\floor{\frac{t_n}{2}}}  \mathds E_{j-t_n}\left[ h_{i,j}(X_i,X_j)\right]-\mathds E_{\pi}\left[h_{i,j}\right] \right|,\\
(2b)&:= \left| \sum_{j=2}^{n} \; \sum_{i=(j-t_n+1)\vee 1}^{j-\floor{\frac{t_n}{2}}}  \mathds E_{\pi}\left[h_{i,j}\right] -\mathds E\left[ h_{i,j}( X_i, X_j)\right]\right|,\\
 \mathrm{and} \quad(2c)&:=\left| \sum_{j=2}^{n} \; \sum_{i=(j-\floor{\frac{t_n}{2}}+1)\vee 1}^{j-1}  \mathds E_{j-t_n}\left[ h_{i,j}(X_i,X_j)\right]-\mathds E\left[h_{i,j}( X_i, X_j)]\right] \right|.
\end{align*}

\raisebox{1.22ex}{\resizebox{!}{1ex}{ \dbend}} The only place where we use the stationarity of the chain is to bound the terms $(2b)$ and $(2c)$ by writing that $\mathds E[h_{i,j}(X_i,X_j)]=\int_{x_i}\int_{x_j}d\pi(x_i)P^{j-i}(x_i,dx_j)$. Both are bounded using similar ideas, that is why we show here how to deal with $(2b)$ and we postpone the proof of Lemma~\ref{lemma:2a-2c} to Section \ref{proof:lemma-2a-2c}. %Appendix~\ref{proof:lemma-2a-2c}.
\begin{lemma}\label{lemma:2a-2c} It holds $\qquad (2a)\leq LAt_n \quad \text{and}\quad  (2c) \leq 2LAt_n^2.$
\end{lemma}
Let us now bound the term $(2b)$,
\begin{align*}
(2b)&:= \left| \sum_{j=2}^{n} \; \sum_{i=(j-t_n+1)\vee 1}^{j-\floor{\frac{t_n}{2}}}  \mathds E_{\pi}\left[h_{i,j}\right] -\mathds E\left[ h_{i,j}( X_i, X_j)\right]\right|\\
&\leq\sum_{j=2}^{n} \; \sum_{i=(j-t_n+1)\vee 1}^{j-\floor{\frac{t_n}{2}}}  \left| \int_{x_i} \int_{x_j}  h_{i,j}(x_i,x_j) d\pi(x_i)\left( d\pi(x_j) -P^{j-i}(x_i,dx_j)\right)\right|\\
&\leq\sum_{j=2}^{n} \; \sum_{i=(j-t_n+1)\vee 1}^{j-\floor{\frac{t_n}{2}}}  \|h_{i,j}\|_{\infty} \underbrace{\int_{x_i}   d\pi(x_i)}_{=1} \underbrace{\sup_{y}\int_{x_j}\left| d\pi(x_j) - P^{j-i}(y,dx_j) \right|}_{= \sup_y\| P^{j-i}(y,\cdot)-\pi\|_{\mathrm{TV}}}\\
%&\leq\sum_{j=2}^{n} \; \sum_{i=(j-t_n+1)\vee 1}^{j-\floor{\frac{t_n}{2}}}  \|h_{i,j}\|_{\infty} \sup_y\| P^{j-i}(y,\cdot)-\pi\|_{\mathrm{TV}}\\
&\leq\sum_{j=2}^{n} \; \sum_{i=(j-t_n+1)\vee 1}^{j-\floor{\frac{t_n}{2}}} \|h_{i,j}\|_{\infty}L\rho^{j-i}\leq\sum_{j=2}^{n} \;\sum_{i=(j-t_n+1)\vee 1}^{j-\floor{\frac{t_n}{2}}} \|h_{i,j}\|_{\infty}L\rho^{t_n/2}\leq LAt_n,
\end{align*}
where we used that~$\rho^{t_n/2}\leq\rho^{r\log(n)/2} = n^{r \log(\rho)/2}\leq n^{-1}~$. Indeed~$1+r\log(\rho)/2<0$ because we choose~$r~$ such that~$r>2(\log(1/\rho))^{-1}.$ Coming back to Eq.\eqref{eq:Rstat}, we deduce that $R^{(t_n)}_{\mathrm{stat}}(n)\leq AL \left( 2+ 2t_n+ 2t_n^2\right)$ which concludes the proof of Proposition~\ref{mainproposition-remaining}.

\section*{Acknowledgements}
This work was supported by a grant from Région Ile de France. The authors gratefully acknowledge the constructive comments and recommendations of the editor and of the anonymous reviewers
which definitely help to improve the readability and quality of the paper. 

\clearpage
\bibliography{sample} 

\clearpage

\appendix

\appendix

%For the sake of completeness, we provide in this supplementary material complete proofs of some Lemmas and we recall important definitions and properties on Markov chains useful for our paper.

% \begin{itemize}
% \item \underline{Section~\ref{appendix}: Definitions and properties on Markov chains}\\
% In the section, we recall useful definitions and results on Markov chains regarding ergodicity, spectral gaps and the splitting method.\smallskip

% %\item \underline{Section~\ref{apdx:han}: Connections with the literature}\\
% %We establish a clear connection between Theorem~\ref{mainthm2} and the exponential inequality from~\cite{shen20} which also provides a concentration inequality for U-statistics in a dependent framework. \smallskip
% \item \underline{Section~\ref{proof-mainthm3}: Proof of Theorem~\ref{mainthm3}}\\
% We provide the small modifications in the proof of Theorem~\ref{mainthm} allowing to get the Bernstein-type concentration result of Theorem~\ref{mainthm3}.\smallskip

% \item \underline{Section~\ref{sec:add-proofs}: Additional proofs}\\
% This section contains the proofs of some Lemmas useful to establish Theorems~\ref{mainthm} and~\ref{mainthm2}. \ref{mainthm3}.
% \end{itemize}

\setlength\parindent{0pt}

\label{appendix}

\section{Proof of Lemma~\ref{lemma:talagrand}}

\label{proof-lemma-talagrand}
In the section, we show that in the proof of Proposition~\ref{mainproposition}, we can use the concentration inequality for the supremum of an empirical process of~\cite[Theorem~3]{samson2000concentration}.

Let us consider the sequence of random variables~$W=(W_1,\dots,W_n)$ on the probability space~$(\Omega,\mathcal A,\mathds P)$ taking values in the measurable space~$ S=E\times F^{n-1}$ where~$F$ is the subset of the set~$\mathcal F(E,\mathds R)$ of all measurable functions from~$(E,\Sigma)$ to~$(\mathds R,\mathcal B(\mathds R))$ that are bounded by~$A$. Note that
\[ \{0_{\mathcal F(E,\mathds R)}\} \cup\{ p_{i,j}(x,\cdot) \; : \; x\in E,\; i,j \in [n]\}\subset F.\]
We define $\mathcal P := \left\{ D \in \mathcal P(F) \; : \; \forall i\in [n-1], \forall j \in \{i+1, \dots, n\}, \; f_{i,j}^{-1}(D) \in \Sigma  \right\}$
where~$\mathcal P(F)$ is the powerset of~$F$ and where~$\forall i\in [n-1], \;  \forall j \in \{i+1, \dots, n\}$, \begin{align*}f_{i,j}:(E,\Sigma) &\to (F,\mathcal P(F))\qquad 
x \mapsto p_{i,j}(x,\cdot).
\end{align*}

Then we have the following straightforward result.
\vspace{-0.25cm}
\begin{lemma}
\label{lemma:measurability}
$\mathcal P$ is a~$\sigma$-algebra on~$F$.
\end{lemma}
\vspace{-0.25cm}

In the following, we endow the space~$F$ with the~$\sigma$-algebra~$\mathcal P$ and we consider on~$ S$ the product~$\sigma$-algebra given by $\mathcal S:=\sigma\left( \left\{C\times D_2 \times \dots \times D_n \; : \; C\in \Sigma, \; D_j \in \mathcal P \;\forall j \in \{2,\dots,n\} \right\}\right) .$\\
For all~$i \in [n]$, we define~$W_i$ by
\[W_i:=\big(X_i,\underbrace{0,\dots,0}_{(i-1) \text{ times}},p_{i,i+1}(X_i,\cdot),p_{i,i+2}(X_i,\cdot), \dots, p_{i,n}(X_i,\cdot)\big).\]
Hence for all~$i \in [n]$,~$W_i$ is~$\sigma(X_i)$-measurable. 
Let us consider for any~$i \in [n-1]$,
\begin{align*}
\Phi_i : (E,\Sigma) &\to (S,\mathcal S) \qquad
\text{such that} \quad \forall x \in E , \; \Phi_i(x)=  \big(x,\underbrace{0_{F},\dots,0_{F}}_{(i-1)\; \text{times}},p_{i,i+1}(x,\cdot),\dots ,p_{i,n}(x,\cdot)\big).
\end{align*}
\vspace{-0.25cm}

Then, one can directly see that for all~$i\in [n-1]$,~$W_i = \Phi_i(X_i)$ and by construction of~$\mathcal P$ and~$\mathcal S$,~$\Phi_i$ is measurable. Indeed, each coordinate of~$\Phi_i$ is measurable by construction of~$\mathcal P$ and this ensures that~$\Phi_i$ is measurable thanks to the following Lemma.

\vspace{-0.25cm}

\begin{lemma} (See~\cite[Lemma 4.49]{aliprantis})
Let~$(X,\Sigma),(X_1,\Sigma_1)$ and~$(X_2,\Sigma_2)$ be measurable spaces, and let~$f_1:X\to X_1$ and ~$f_2:X\to X_2$. Define ~$f:X\to X_1\times X_2$ by~$f(x)=(f_1(x),f_2(x))$. Then ~$f:(X,\Sigma)\to (X_1\times X_2,\Sigma_1\otimes \Sigma_2)$ is measurable if and only if the two functions~$f_1:(X,\Sigma)\to (X_1,\Sigma_1)$ and ~$f_2:(X,\Sigma)\to(X_2,\Sigma_2)$ are both measurable.
\end{lemma}
Then it holds for any~$i \in \{2,\dots,n-1\}$ and any~$G\in \mathcal S$,
\begin{align}
&\mathds P\left( W_{i} \in G\; | \; W_{i-1}   \right)=\mathds P\left( \Phi_{i}(X_i) \in G \; | \; W_{i-1}   \right)=\mathds P\left( \Phi_{i}(X_i)  \in G \; | \; X_{i-1}   \right)\notag\\
&=\mathds P\left( X_i  \in \Phi_i^{-1}\left( G\right) \; | \; X_{i-1}   \right)= P\left(X_{i-1},\Phi_i^{-1}\left( G\right)\right)= \left[(\Phi_i)_{\#}P(X_{i-1},\cdot)\right] \; \left( G \right), \label{pushforward}
\end{align}
where~$(\Phi_i)_{\#}P(X_{i-1},\cdot)$ denotes the pushforward measure of the measure~$P(X_{i-1},\cdot)$ by the measurable map~$\Phi_i$. We deduce that~$W_i$ is non-homogeneous Markov chain. Moreover,~\eqref{pushforward} proves that the transition kernel of the Markov chain~$(W_k)_k$ from state~$i-1$ to state~$i$ is given by~$K^{(i-1,i)}$ where for all~$ (x,p_2,\dots,p_n) \in S$ and for all~$G \in \mathcal S$,
\begin{align*}
  &K^{(i-1,i)}((x,p_2,\dots,p_n),G) =\left[ (\Phi_i)_{\#}P(x,\cdot)\right]\left(G\right).
\end{align*}
One can easily generalize this notation. Let us consider some~$i,j\in[n]$ with~$i<j$ and let us denote by~$K^{(i,j)}$ the transition kernel of the Markov chain~$(W_k)_k$ from state~$i$ to state~$j$. Then for all~$x \in E$, for all~$p_2,\dots,p_n \in F$ and for all~$ G\in \mathcal S,~$
\begin{align*}
&K^{(i,j)}( (x,p_2,\dots,p_n)\;,G)= \left[ (\Phi_j)_{\#}P^{j-i}(x,\cdot)\right]\left(G\right),  \end{align*}
We introduce the mixing matrix~$\Gamma=(\gamma_{i,j})_{1\leq i,j\leq n-1}$ where coefficients are defined by
\[\gamma_{i,j} := \sup_{w_i\in  S} \sup_{z_i\in  S} \|\mathcal L(W_j|W_i=w_i)-\mathcal L(W_j|W_i=z_i)\|_{\mathrm{TV}}. \]
For any~$w \in  S=E\times F^{n-1}$, we denote by~$w^{(1)}$ the first coordinate of the vector~$w$. Hence,~$w^{(1)}$ is an element of~$E$. Then
\begin{align*}
\gamma_{i,j} &= \sup_{w_i \in  S} \sup_{z_i \in  S} \sup_{G \in \mathcal S} \left| \left[ (\Phi_j)_{\#}P^{j-i}(w_i^{(1)},\cdot)\right]\left(G\right) - \left[ (\Phi_j)_{\#}P^{j-i}(z_i^{(1)},\cdot)\right]\left(G\right) \right|\\
 &=\sup_{w_i \in  S} \sup_{z_i \in  S} \sup_{G \in \mathcal S} \left|  P^{j-i}\left(w_i^{(1)},\Phi_j^{-1}\left(G\right) \right)-  P^{j-i}\left(z_i^{(1)},\Phi_j^{-1}\left(G\right) \right) \right|\\
 &\leq \sup_{w_i \in  S} \sup_{z_i \in  S} \sup_{C \in \Sigma}  \left|  P^{j-i}\left(w_i^{(1)},C \right)-  P^{j-i}\left(z_i^{(1)},C \right) \right|\\
 &=\sup_{x_i \in E} \sup_{x'_i \in E} \sup_{C \in \Sigma}  \left|  P^{j-i}\left(x_i,C \right)-  P^{j-i}\left(x'_i,C \right) \right|=\sup_{x_i \in E} \sup_{x'_i\in E} \|P^{j-i}(x_i,\cdot) - \pi(\cdot)+\pi(\cdot)- P^{j-i}(x'_i,\cdot)\|_{\mathrm{TV}}\\
&\leq \sup_{x_i \in E}  \|P^{j-i}(x_i,\cdot) - \pi(\cdot)\|_{\mathrm{TV}} + \sup_{x'_i\in E} \|P^{j-i}(x'_i,\cdot) - \pi(\cdot)\|_{\mathrm{TV}}\leq 2L \rho^{j-i},
\end{align*}
where in the first inequality we used that~$\Phi_j:(E,\Sigma) \to (S,\mathcal S)$ is measurable and in the last inequality we used the uniform ergodicity of the Markov chain~$(X_i)_{i\geq 1}.$ We deduce that 
\[\|\Gamma\| \leq 2L \left\| \mathrm{Id} + \sum_{l=1}^{n-1}\rho^lN_l\right\|,\quad N_l=\left(n^{(l)}_{i,j}\right)_{1\leq i,j \leq n-1} \; \text{with} \;n^{(l)}_{i,j}= \left\{
    \begin{array}{ll}
        1 & \mbox{if } j-i=l\\
        0 & \mbox{otherwise.}
    \end{array}
\right.. \]
Note that $N_l$ is a nilpotent matrix of order~$l$. Since for each~$1\leq l \leq n-1$,~$\|N_l\|\leq1$, it follows from the triangular inequality that
\[\|\Gamma\| \leq 2L \sum_{l=0}^{n-1} \rho^l\leq \frac{2L}{1-\rho}.\]
To conclude the proof and get the concentration result stated in Lemma~\ref{lemma:talagrand}, one only needs to apply~\cite[Theorem~3]{samson2000concentration} with the class of functions~$\mathcal F$ and with the Markov chain~$(W_k)_k$. Let us recall that~$\mathcal F$ is defined by ~$\mathcal{F}=\{f_{\xi} : \sum_{j=2}^n \mathds E|\xi_j(X'_j)|^{k/(k-1)}=1\}$ where for any~$\xi=(\xi_2, \dots, \xi_n) \in \prod_{i=2}^n L^{k/(k-1)}(\nu)$, \begin{equation*}\forall w=(x,p_2,\dots,p_n) \in E\times F^{n-1}, \quad f_{\xi}(w) = \sum_{j=2}^n \int p_j(y) \xi_j(y) d\nu(y).\end{equation*}

\section{Proof of Lemma~\ref{lemma:bound-Vnk-proba}}
\label{subsec:proof-bound-Vnk-proba}

\textbf{Bounding~$b_k$.} Using Hölder's inequality we have,
\begin{align*}
b_k&=\sup_{w \in S} \sup_{f_{\xi} \in \mathcal F} |f_{\xi}(w)|= \sup_{(p_2,\dots,p_n) \in F^{n-1}} \sup_{\xi \in \mathcal L_k} \sum_{j=2}^n \mathds E[p_j(X'_j)\xi_j(X'_j)] \\
&\leq    \sup_{(p_2,\dots,p_n) \in F^{n-1}} \underset{\sum_{j=2}^n \mathds E |\xi_j(X'_j)|^{k/(k-1)}=1}{\sup}  \sum_{j=2}^n \left(\mathds E\left|p_{j}(X'_j)\right|^{k}\right)^{1/k}\left(\mathds E\left|\xi_j(X'_j)\right|^{k/(k-1)} \right)^{(k-1)/k}\\
&\leq \sup_{(p_2,\dots,p_n) \in F^{n-1}} \underset{\sum_{j=2}^n \mathds E |\xi_j(X'_j)|^{k/(k-1)}=1}{\sup}   \left( \sum_{j=2}^n\mathds E\left|p_{j}(X'_j)\right|^{k}\right)^{1/k}\left( \sum_{j=2}^n\mathds E\left|\xi_j(X'_j)\right|^{k/(k-1)} \right)^{(k-1)/k}\\
&\leq  \sup_{(p_2,\dots,p_n) \in F^{n-1}}  \left( \sum_{j=2}^n\mathds E\left|p_{j}(X'_j)\right|^{k}\right)^{1/k}\leq ((nA^2)A^{k-2})^{1/k},
\end{align*}
where~$A:=  2 \max_{i,j}\|h_{i,j}\|_{\infty}$ which satisfies~$\max_{i,j}\|p_{i,j}\|_{\infty}\leq A.$ Here, we used that~$F$ is the set of measurable functions from~$(E,\Sigma)$ to~$(\mathds R, \mathcal B(\mathds R))$ bounded by~$A$.

%\[B^2:=\max \left\{ \max_j\left\| \sum_{i=1}^{j-1} \mathds E_{i-1} h_{i,j}^2(X_i,x) \right\|_{\infty} \quad ,   \max_i\left\| \sum_{j=i+1}^{n} \mathds E_{j-1} h_{i,j}^2(x,X_j) \right\|_{\infty}   \right\}\]

\textbf{Bounding the variance.}
\vspace{-0.4cm}

\begin{align*}
\sigma^2_k&= \mathds E\left[ \sum_{i=1}^{n-1}\underset{f_{\xi} \in \mathcal{F}}{\sup} f_{\xi}(W_i)^2 \right]= \sum_{i=1}^{n-1} \mathds E\left[   \underset{\xi \in \mathcal L_k}{\sup} \left( \sum_{j=i+1}^n \mathds E_{X'_j}\left[ p_{i,j}(X_i,X'_j)\xi_j(X'_j)\right]\right)^2 \right] \\
& = \sum_{i=1}^{n-1} \mathds E\left[    \left(\underset{\xi \in \mathcal L_k}{\sup} \left| \sum_{j=i+1}^n \mathds E_{X'_j}\left[ p_{i,j}(X_i,X'_j)\xi_j(X'_j)\right]\right|\right)^2 \right] \leq n\left(\mathfrak B_0^2A^{k-2}\right)^{2/k} ,
\end{align*}

\noindent
where the last inequality comes from the following (where we use twice Holder's inequality),
\vspace{-0.4cm}

\begin{align*}
&\underset{\xi \in \mathcal L_k}{\sup} \left| \sum_{j=i+1}^n \mathds E_{X'_j}\left[ p_{i,j}(X_i,X'_j)\xi_j(X'_j)\right]\right| \\
\leq\quad &  \underset{\xi \in \mathcal{L}_k}{\sup} \sum_{j=i+1}^n \left(\mathds E_{X'_j}\left|p_{i,j}(X_i,X'_j)\right|^{k}\right)^{1/k}\left(\mathds E\left|\xi_j(X'_j)\right|^{k/(k-1)} \right)^{(k-1)/k}\\
\leq\quad &   \underset{\sum_{j=2}^n \mathds E |\xi_j(X'_j)|^{k/(k-1)}=1}{\sup}   \left( \sum_{j=i+1}^n\mathds E_{X'_j}\left|p_{i,j}(X_i,X'_j)\right|^{k}\right)^{1/k}\left( \sum_{j=i+1}^n\mathds E_{X'_j}\left|\xi_j(X'_j)\right|^{k/(k-1)} \right)^{(k-1)/k}\\
\leq\quad &   \left( \sum_{j=i+1}^n\mathds E_{X'_j}\left|p_{i,j}(X_i,X'_j)\right|^{k}\right)^{1/k}\leq\left( \mathfrak B_0^2 A^{k-2}\right)^{1/k},\quad \text{where }\mathfrak B_0\text{ is defined in }\eqref{eq:B0}. 
\end{align*}

\noindent
Using Lemma~\ref{lemma:inequality} twice and the bounds obtained on~$b_k$ and~$\sigma_k^2$ gives for~$u>0$,
\begingroup
 \allowdisplaybreaks
\begin{align*}
&\left[(2\delta_M)^{1/k} \mathds E[Z]+(2\delta_M)^{1/k} \sigma_k 3\|\Gamma\|\sqrt{ku}+(2\delta_M)^{1/k} k8\|\Gamma\|^2b_ku\right]^{k} \\
\leq \quad & \Bigg[ (2\delta_M)^{1/k}\mathds E[Z]+(2\delta_M)^{1/k}3\|\Gamma\| ( \mathfrak B_0^2A^{k-2})^{1/k} \sqrt{n ku}+(2\delta_M)^{1/k}8\|\Gamma\|^2 ((nA^2)A^{k-2})^{1/k}ku\Bigg]^{k}  \\
\leq \quad &(1+\epsilon)^{k-1}2\delta_M \left(\mathds E[Z]\right)^k +(1+\epsilon^{-1})^{k-1}\Bigg[(2\delta_M)^{1/k} 8\|\Gamma\|^2( (nA^2)A^{k-2})^{1/k}ku \\
&\quad+(2\delta_M)^{1/k} 3\|\Gamma\|(  \mathfrak B_0^2A^{k-2})^{1/k}\sqrt{n ku}\Bigg]^{k} \\
\leq \quad & (1+\epsilon)^{k-1}2\delta_M \left(\mathds E[Z]\right)^k +2\delta_M(1+\epsilon^{-1})^{2k-2}\left(8\|\Gamma\|^2\right)^k(nA^2)A^{k-2}( ku)^k\\
& \quad +(1+\epsilon)^{k-1}(1+\epsilon^{-1})^{k-1}2\delta_M \left(3\|\Gamma\|\right)^k  \mathfrak B_0^2A^{k-2} (n ku)^{k/2}.
\end{align*}
\endgroup

\section{Bounding $(\mathds E [Z])^k$ under Assumption~\ref{assumption0}$.(ii)$}
\label{apdx:EZ-chi}

In this section, we only provide the part of the proof of Proposition~\ref{mainproposition} that needs to be modified to get the result when the kernels~$h_{i,j}$ depend on both~$i$ and~$j$ and when Assumption~\ref{assumption0}$.(ii)$ holds. Keeping the notations of the proof of Proposition~\ref{mainproposition}, we only want to bound~$(\mathds E[Z])^k $ (and thus $a_1$) using a different concentration result that can allow to deal with kernel functions~$h_{i,j}$ that might depend on~$i$. We will use Proposition~\ref{bernstein-markov}.

\begin{proposition}  \label{bernstein-markov} (cf.~\cite[Theorem 1]{jiang2018bernsteins})
Suppose that the sequence~$(X_i)_{i\geq1}$ is a Markov chain satisfying Assumptions~\ref{assumption1} and~\ref{assumption0}$.(ii)$ with invariant distribution~$\pi$ and with an absolute spectral gap~$1-\lambda>0$. Let us consider some~$n \in \mathds N^*$ and bounded real valued functions~$(f_i)_{1\leq i \leq n}$ such that for any~$i \in \{1, \dots,n\}$,~$\int f_i(x) d\pi(x)=0$ and~$\|f_i\|_{\infty}\leq c$ for some~$c>0$. Let~$\sigma^2 = \sum_{i=1}^n \int f_i^2(x) d\pi(x) /n$.
Then for any~$0 \leq t < (1-\lambda)/(5cq)$,
\begin{equation}
\mathds E_{\chi}\left[ e^{t\sum_{i=1}^n f_i(X_i)} \right] \leq \left\| \frac{d\chi}{d\pi} \right\|_{\pi,p}\exp\left( \frac{n \sigma^2}{qc^2}(e^{tqc}-tqc-1)+\frac{n\sigma^2\lambda q t^2}{1-\lambda-5cqt}  \right), \label{majo-laplace}
\end{equation} where~$q$ is the constant introduced in Assumption~\ref{assumption0}$.(ii)$. Moreover for any~$u \geq 0$ it holds
%\[\mathds P \left( \sum_{i=1}^n f_i(X_i) \geq \epsilon \right) \leq  \left\| \frac{d\chi}{d\pi} \right\|_{\pi,p} \exp\left(- \frac{\epsilon^2/(2q)}{A_2n\sigma^2+A_1c\epsilon}\right) ,\]
\[\mathds P \left( \frac{1}{n} \sum_{i=1}^n f_i(X_i)> \frac{2quA_1 c}{n} + \sqrt{\frac{2quA_2\sigma^2}{n}} \right) \leq  \left\| \frac{d\chi}{d\pi} \right\|_{\pi,p} e^{-u}.\]
where~$A_2 := \frac{1+\lambda}{1-\lambda}$ and~$A_1:=\frac13 \mathds1_{\lambda=0}+ \frac{5}{1-\lambda}\mathds 1_{\lambda >0}$.

\end{proposition}

Note that the Bernstein inequality~\cite[Theorem 1]{jiang2018bernsteins} is given for stationary chains but one can easily extend this result to obtain Proposition~\ref{bernstein-markov} by working under the milder assumption Assumption~\ref{assumption0}$.(ii)$ following for example the approach used in~\cite[Theorem 2.3]{FJS18}. Let us recall that
\begin{align*}
(\mathds E[Z])^k &\leq  \mathds E[Z^k] \quad \text{(Using Jensen's inequality)}\\
&= \mathds E\left[\left(  \underset{f_{\xi} \in \mathcal{F}}{\sup}  \sum_{i=1}^{n-1} f_{\xi}(X_i) \right)^k\right] = \mathds E\left[\left(  \underset{f_{\xi} \in \mathcal{F}}{\sup} \sum_{i=1}^{n-1} \sum_{j=i+1}^{n} \mathds E_{j-1}[p_{i,j}(X_i,X'_j)\xi_j(X'_j)] \right)^k\right] \\
&= \mathds E\left[\sum_{j=2}^{n} \mathds E_{j-1} \Bigg|\sum_{i=1}^{j-1} p_{i,j}(X_i,X'_j)\Bigg|^k\right] \quad \text{(Using Lemma~\ref{duality})}\\
&= \mathds E\left[\sum_{j=2}^{n} \mathds E_{|X'} \Bigg|\sum_{i=1}^{j-1} p_{i,j}(X_i,X'_j)\Bigg|^k\right] .\end{align*}
 Thus we have
 \[a_1= \frac{2\delta_M}{1+\epsilon} \mathds{E}\sum_{j=2}^n\Bigg( \mathds{E}_{|X'} \left[ e^{ \alpha (1+\epsilon)K |C^{(j)}|}\right] -\alpha (1+\epsilon)K\mathds E_{|X'}\left[|C^{(j)}|\right] -1\Bigg),
\]
where~$C^{(j)} = \sum_{i=1}^{j-1} p_{i,j}(X_i,X'_j)$ and where the notation~$\mathds{E}_{|X'}$ refers to the expectation conditionally to the~$\sigma$-algebra~$\sigma(X'_{2},\dots,X'_n)$.

Now we use a symmetrization trick: since~$e^x-x-1\geq0$ for all~$x$ and since~$e^{a|x|}+e^{-a|x|} = e^{ax}+e^{-ax}$, adding~$\mathds{E}_{|X'} [ \exp\left( -\alpha (1+\epsilon)K|C^{(j)}|\right)] +\alpha (1+\epsilon)K\mathds E_{|X'} [|C^{(j)}|] -1$ to~$a_1$ gives
\begin{equation}a_1 \leq \frac{2\delta_M}{1+\epsilon} \mathds{E}\sum_{j=2}^n \Bigg( \mathds{E}_{|X'} [ e^{ \alpha (1+\epsilon)K C^{(j)}}]  -1+\mathds{E}_{|X'} [e^{ -\alpha (1+\epsilon) KC^{(j)}}]  -1\Bigg).\label{majo-a1-bern}\end{equation}

Let us consider some~$j \in \{2, \dots,n\}$. Conditionally on~$\sigma(X'_2,\dots,X'_n)$, ~$ C^{(j)}$ is a sum of bounded functions (by~$A$) depending on the Markov chain. We denote \[v_j(X'_j) =\sum_{i=1}^{j-1} \mathds{E}_{X_i \sim \pi}[p_{i,j}^2(X_i,X'_j) | X'_j] \leq \mathfrak B_0^2\] and~$V=\sum_{j=2}^n \mathds E v_j^k(X'_j)\leq  C_n^2 \mathfrak B_0^{2(k-1)}$ (with~$C_n^2=\sum_{j=2}^n \sum_{i=1}^{j-1} \mathds E [ p_{i,j}^2(X_i,X'_j)]$).

\noindent
Remark that
\begin{align*}
\mathds{E}_{X_i \sim \pi}[p_{i,j}(X_i,X'_j) | X'_j]&=\mathds{E}_{X_i \sim \pi}\left[h_{i,j}(X_i,X'_j) -\mathds E_{\tilde X\sim \pi}[h_{i,j}(X_i,\tilde X)] \;|\; X'_j \right]\\
&= \int_{x' } \left( \int_{x_i }(h_{i,j}(x_i,X'_j)- h_{i,j}(x_i,\tilde x))d\pi(x_i)\right)d\pi(\tilde x)=0,
\end{align*}
where the last equality comes from Assumption~\ref{assumption4}. We use the Bernstein inequality for Markov chain (see Proposition~\ref{bernstein-markov}). Notice from Taylor expansion that~$(1-p/3)(e^p-p-1) \leq p^2/2$ for all~$p\geq 0$. Applying~\eqref{majo-laplace} with~$t=\alpha(1+\epsilon)K$ and~$c=A$, we get that for~$\alpha < [(1+\epsilon)K\sqrt{q}(A\sqrt{q}/3+\mathfrak B_0\sqrt{3/2})]^{-1}\wedge [ (1-\lambda)^{-1/2}(1+\epsilon) K \sqrt{q} \left( 5A\sqrt{q}(1-\lambda)^{-1/2}+\sqrt{3\lambda}\mathfrak B_0  \right) ]^{-1}$,
\begin{align*}& \mathds{E}_{|X'} [ e^{ \alpha (1+\epsilon)K |C^{(j)}|}]\leq   2\left\| \frac{d\chi}{d\pi} \right\|_{\pi,p} \times   \mathds E_{|X'} \left[ \exp\left( \frac{\alpha^2(1+\epsilon)^2K^2qv_j(X'_j)}{2-2Aq\alpha(1+\epsilon)K/3}+\frac{v_j(X'_j)\lambda \alpha^2(1+\epsilon)^2K^2q}{1-\lambda-5 \alpha(1+\epsilon)KAq}\right)   \right].\end{align*}

\noindent
Considering~$\alpha < [(1+\epsilon)K\sqrt{q}(A\sqrt{q}/3+\mathfrak B_0\sqrt{3/2})]^{-1}\wedge [(1-\lambda)^{-1/2}(1+\epsilon) K \sqrt{q} \left( 5A\sqrt{q}(1-\lambda)^{-1/2}+\sqrt{3\lambda}\mathfrak B_0  \right) ]^{-1}$,~$\epsilon<1$ and using~\eqref{majo-a1-bern}, this leads to
\begin{align*}
&\frac{a_1}{2\left\| \frac{d\chi}{d\pi} \right\|_{\pi,p}} \leq   \frac{2\delta_M}{1+\epsilon} \sum_{j=2}^n \mathds E \left[ \exp\left( \frac{\alpha^2(1+\epsilon)^2K^2qv_j(X'_j)}{2-2Aq\alpha(1+\epsilon)K/3}+\frac{v_j(X'_j)\lambda \alpha^2(1+\epsilon)^2K^2q}{1-\lambda-5 \alpha(1+\epsilon)KAq}\right)-1   \right]\\
 & = \frac{2\delta_M}{1+\epsilon} \sum_{j=2}^n \sum_{k=1}^{\infty} \frac{1}{k!} \left( \frac{\alpha^2(1+\epsilon)^2K^2qv_j(X'_j)}{2-2Aq\alpha(1+\epsilon)K/3}+\frac{v_j(X'_j)\lambda \alpha^2(1+\epsilon)^2K^2q}{1-\lambda-5 \alpha(1+\epsilon)KAq}\right)^k \\
  & = \frac{2\delta_M}{1+\epsilon} \sum_{j=2}^n \sum_{k=1}^{\infty} \frac{1}{k!}\left(\frac{3}{2}\right)^{k-1} \left( \frac{\alpha^2(1+\epsilon)^2K^2qv_j(X'_j)}{2-2Aq\alpha(1+\epsilon)K/3}\right)^k\\
  &\qquad +\frac{2\delta_M}{1+\epsilon} \sum_{j=2}^n \sum_{k=1}^{\infty} \frac{1}{k!}3^{k-1} \left( \frac{v_j(X'_j)\lambda \alpha^2(1+\epsilon)^2K^2q}{1-\lambda-5 \alpha(1+\epsilon)KAq}\right)^k \quad \text{(Using Lemma~\ref{lemma:inequality})} \\
 & \leq  \frac{\delta_M}{3(1+\epsilon)} \sum_{k=1}^{\infty}  \frac{3^k\alpha^{2k}(1+\epsilon)^{2k}K^{2k}q^kV}{(4-4Aq\alpha(1+\epsilon)K/3)^k}+\frac{2\delta_M}{3(1+\epsilon)} \sum_{k=1}^{\infty}  \frac{3^kV\lambda^k \alpha^{2k}(1+\epsilon)^{2k}K^{2k}q^k}{\left(1-\lambda-5 \alpha(1+\epsilon)KAq\right)^k} \\
  & \leq    \frac{\delta_M}{3(1+\epsilon)} \sum_{k=1}^{\infty}  \frac{3^k\alpha^{2k}(1+\epsilon)^{2k}K^{2k}q^kC_n^2\mathfrak B_0^{2(k-1)}}{(2-2Aq\alpha(1+\epsilon)K/3)^k}+\frac{2\delta_M}{3(1+\epsilon)} \sum_{k=1}^{\infty}  \frac{3^kC_n^2\mathfrak B_0^{2(k-1)}\lambda^k \alpha^{2k}(1+\epsilon)^{2k}K^{2k}q^k}{\left(1-\lambda-5 \alpha(1+\epsilon)KAq\right)^k}\\
 &=   \frac{(1+\epsilon)C_n^2\alpha^2K^2\delta_Mq}{2-2Aq\alpha (1+\epsilon)K/3-3\alpha^2(1+\epsilon)^2K^2\mathfrak B_0^2q}+\frac{2\delta_M C_n^2 \lambda \alpha^2 (1+\epsilon)K^2q}{1-\lambda-5\alpha (1+\epsilon)KAq-3\mathfrak B_0^2\lambda \alpha^2 (1+\epsilon)^2K^2q}\\
  &=   \frac{(1+\epsilon)C_n^2\alpha^2K^2\delta_M q/2}{1-Aq\alpha (1+\epsilon)K/3-3\alpha^2(1+\epsilon)^2K^2\mathfrak B_0^2q/2} +\frac{2\delta_M C_n^2 \lambda \alpha^2 (1+\epsilon)K^2q(1-\lambda)^{-1}}{1-5(1-\lambda)^{-1}\alpha  (1+\epsilon)KAq-3\mathfrak B_0^2\lambda (1-\lambda)^{-1} \alpha^2 (1+\epsilon)^2K^2q}\\
 &\leq  \frac{(1+\epsilon)C_n^2\alpha^2K^2\delta_M q/2}{1-\alpha (1+\epsilon)K\sqrt{q}(A\sqrt{q}/3+\mathfrak B_0\sqrt{3/2})}+\frac{2\delta_M C_n^2 \lambda \alpha^2 (1+\epsilon)K^2q(1-\lambda)^{-1}}{1-\alpha(1-\lambda)^{-1/2}(1+\epsilon) K \sqrt{q} \left( 5A\sqrt{q}(1-\lambda)^{-1/2}+\sqrt{3\lambda}\mathfrak B_0  \right)}.
\end{align*}

From this bound on~$a_1$, one can follow the steps of the proof of Proposition~\ref{mainproposition} to conclude.

\section{Proof of Lemma~\ref{lemma:bound-a1a2a3}}
\label{apdx:bound-a1a2a3}

\textbf{Bounding~$a_3$.} Using the inequality~$k! \geq (k/e)^k$, we have,
\begin{align*}a_3 &\leq  2\delta_M\sum_{k\geq2}\alpha^k(1+\epsilon^{-1})^{2k-2}\left(8\|\Gamma\|^2\right)^k  (nA^2)A^{k-2}( eu)^k\\
&= 2 \delta_M\alpha^2\left[\sqrt{n}A(1+\epsilon^{-1})8\|\Gamma\|^2eu\right]^2\sum_{k\geq2}\alpha^{k-2}(1+\epsilon^{-1})^{2(k-2)}\left(8\|\Gamma\|^2\right)^{k-2}  A^{k-2}( eu)^{k-2}\\
&= \frac{2 \delta_M\alpha^2\left[\sqrt{n}A(1+\epsilon^{-1})8\|\Gamma\|^2eu\right]^2}{1-\alpha(1+\epsilon^{-1})^{2}\left(8\|\Gamma\|^2\right)  A eu},\quad \text{for }\alpha < ((1+\epsilon^{-1})^{2}\left(8\|\Gamma\|^2\right)  A eu)^{-1}\,.
\end{align*}

\pagebreak[3]

\textbf{Bounding~$a_2$.} We use the inequality~$k! \geq k^{k/2}$ because~$(k/e)^k >k^{k/2}$ for~$k \geq e^2$ and for~$k$ smaller, the inequality follows by direct verification. Hence,
\begin{align*}
a_2&\leq  2\delta_M \sum_{k\geq2}\alpha^k(2+\epsilon+\epsilon^{-1})^{k-1}\left(3\|\Gamma\|\right)^k    \mathfrak B_0^2A^{k-2} (n u)^{k/2}\\
&=   2\delta_M(2+\epsilon+\epsilon^{-1})\alpha^2\left[3\|\Gamma\|  \mathfrak B_0 \sqrt{nu}\right]^2\sum_{k\geq2}\alpha^{k-2}(2+\epsilon+\epsilon^{-1})^{k-2}\left(3\|\Gamma\|\right)^{k-2}  A^{k-2} (n u)^{(k-2)/2}\\
&= \frac{2\delta_M(2+\epsilon+\epsilon^{-1})\alpha^2\left[3\|\Gamma\|  \mathfrak B_0 \sqrt{nu}\right]^2}{1-\alpha(2+\epsilon+\epsilon^{-1})\left(3\|\Gamma\|\right)  A (n u)^{1/2}}, \quad \text{for}~\alpha < ((2+\epsilon+\epsilon^{-1})\left(3\|\Gamma\|\right)  A (n u)^{1/2})^{-1}.
\end{align*}

\textbf{Bounding~$a_1$.}  Using the bound previously obtained for $(\mathds E[Z])^k$ we get,
\begin{align*}
a_1&= 2\delta_M \sum_{k\geq2}\frac{\alpha^k}{k!}(1+\epsilon)^{k-1} (\mathds E[Z])^k\leq  32\delta_M n\sum_{k\geq2}\alpha^k(1+\epsilon)^{k-1} (KAm\tau)^k\\
&\leq 32\delta_M n \alpha^2 (1+\epsilon)[KAm\tau]^2\sum_{k\geq2}\alpha^{k-2}(1+\epsilon)^{k-2}(KAm \tau)^{k-2}\\
&\leq   \frac{ 32\delta_M n \alpha^2 (1+\epsilon)[KAm\tau]^2}{1-\alpha (1+\epsilon)KAm\tau}  , \quad \text{for }0<\alpha < ( (1+\epsilon)KAm\tau)^{-1}.
\end{align*}

\section{Proof of Lemma~\ref{lemma:induction}}
\label{apdx:proof-lemma-induction}

We prove a concentration result for the term
\[U_{n-1}^{(1)}=\sum_{j=2}^{n-1} \sum_{i=1}^{j-1} h^{(1)}_{i,j}(X_i,X_{j-1},X_{j}),\]
where $h^{(1)}_{i,j}(x,y,z) = \int_w h_{i,j}(x,w)P(z,dw) -\int_w h_{i,j}(x,w)P^2(y,dw),$ using an approach similar to the one from Section~\ref{sec:firstterm}.

\begin{itemize}
\item \underline{Martingale structure} \\
Denoting~$Y_{j}^{(1)} = \sum_{i=1}^{j-1}  h^{(1)}_{i,j}(X_i,X_{j-1},X_{j})$, we have~$U_{n-1}^{(1)}=\sum_{j=2}^{n-1} Y_{j}^{(1)}$ which shows that~$(U_{n}^{(1)})_n$ is a martingale with respect to the~$\sigma$-algebras~$(G_l)_l$. Indeed, we have~$\mathds E_{j-1}[Y_{j}^{(1)}]=0$.

\item \underline{Talagrand's inequality}  To upper-bound~$(V_n^k)_n$, we split it as previously namely
\begin{align*}
\hspace{-0.4cm}V_n^k &:= \sum_{j=2}^{n-1}  \mathds E_{j-1} \left| \sum_{i=1}^{j-1} h^{(1)}_{i,j}(X_i,X_{j-1},X_j)\right|^k \\
&= \sum_{j=2}^{n-1} \mathds E_{j-1} \Bigg| \sum_{i=1}^{j-1} \Bigg(   I^{(1)}_{i,j}(X_i,X_j) -\mathds E_{j-1}[ I^{(1)}_{i,j}(X_i,X_j)] \Bigg)\Bigg|^k,
\end{align*}
where $I^{(1)}_{i,j}(x,z) =\int_w h_{i,j}(x,w)P(z,dw) .$
Using as previously Lemma~\ref{lemma:inequality} with~$\epsilon=1/2$, we get
\begin{align*}
V_n^k &= \sum_{j=2}^{n-1}  \mathds E_{j-1} \Bigg| \sum_{i=1}^{j-1} \Bigg(   I^{(1)}_{i,j}(X_i,X_j) -\mathds E_{\tilde X\sim \pi}[I^{(1)}_{i,j}(X_i,\tilde X)] \\
&\qquad \qquad+\mathds E_{\tilde X\sim \pi}[I^{(1)}_{i,j}(X_i,\tilde X)] - \mathds E_{j-1}[ I^{(1)}_{i,j}(X_i,X_j)] \Bigg)\Bigg|^k\\
& \leq (3/2)^{k-1} \sum_{j=2}^{n-1}  \mathds E_{j-1} \Bigg| \sum_{i=1}^{j-1}\left(  I^{(1)}_{i,j}(X_i,X_j) - \mathds E_{\tilde X\sim \pi}[I^{(1)}_{i,j}(X_i,\tilde X)]\right)\Bigg|^k \\
& \quad +3^{k-1}\sum_{j=2}^{n-1}  \mathds E_{j-1} \Bigg| \sum_{i=1}^{j-1}\left(\mathds E_{\tilde X\sim \pi}[I^{(1)}_{i,j}(X_i,\tilde X)]- \mathds E_{j-1}[ I^{(1)}_{i,j}(X_i,X_j)]\right) \Bigg|^k.
\end{align*}

Again, basic computations and Jensen's inequality lead to \begin{align*}&\sum_{j=2}^{n-1}  \mathds E_{j-1} \Bigg| \sum_{i=1}^{j-1}\left(\mathds E_{\tilde X\sim \pi}[I^{(1)}_{i,j}(X_i,\tilde X)]- \mathds E_{j-1}[ I^{(1)}_{i,j}(X_i,X_j)]\right) \Bigg|^k \\
&=  \sum_{j=2}^{n-1}  \mathds E_{j-1} \Bigg| \sum_{i=1}^{j-1} p^{(1)}_{i,j}(X_i,X_j)\Bigg|^k,\end{align*}
where $p^{(1)}_{i,j}(x,z):=I^{(1)}_{i,j}(x,z) -\mathds E_{\tilde X\sim \pi}[I^{(1)}_{i,j}(x,\tilde X)].$
Hence, using Assumption~\ref{assumption1} and Lemma~\ref{lemma:decoupling} exactly like in the previous section, we get (for $(X'_j)_j$ i.i.d. with distribution $\nu$)
\begin{align*}
V_n^k &= 2\times 3^{k-1}\sum_{j=2}^{n-1}  \mathds E_{j-1} \Bigg| \sum_{i=1}^{j-1} p^{(1)}_{i,j}(X_i,X_j)\Bigg|^k\leq 2\times 3^{k-1}\delta_M\sum_{j=2}^{n-1}  \mathds E_{X'_j} \Bigg| \sum_{i=1}^{j-1} p^{(1)}_{i,j}(X_i,X_j')\Bigg|^k.
\end{align*}

Then, one can use the same duality trick to show that the~$V_n^k$ can be controlled using the supremum of a sum of functions of the Markov chain~$(X_i)_{i\geq1}$ using~\cite[Theorem~3]{samson2000concentration}.

\item \underline{Bounding~$\exp(w_n^k \alpha^k / k!)$}\\
The terms~$a_2$ and~$a_3$ (see Eq.\eqref{eq:defai}) can be bounded in a similar way. For the term~$a_1$, we only need to show that~$p^{(1)}_{i,j}$ satisfies~$\mathds E_{X_i \sim \pi}|p^{(1)}_{i,j}(X_i,z)]=0, \; \forall z\in E$ in order to apply as previously a Bernstein's type inequality.
\begin{align*}
\mathds E_{X_i \sim \pi}|p^{(1)}_{i,j}(X_i,z)]&=  \int_{x_i} d\pi(x_i)  \int_w h_{i,j}(x_i,w)P(z,dw)- \mathds E_{X\sim \pi}\mathds E_{\tilde X\sim \pi}[I^{(1)}_{i,j}(X,\tilde X)]\\
&= \mathds E_{\pi}[h_{i,j}] - \mathds E_{\pi}[h_{i,j}] \quad \text{(Using Assumption~\ref{assumption4})}\\
&= 0.
\end{align*}

\item \underline{Conclusion of the proof}\\
Let us consider the quantities~$A_1$ and~$\mathfrak B_1$ defined as the counterparts of~$A$ and~$\mathfrak B_0$ (see~\eqref{eq:B0}) by replacing the functions~$\left(p_{i,j}\right)_{i,j}$ by~$\left(p^{(1)}_{i,j}\right)_{i,j}$. One can easily see that~$A_1=A$. Let us give details about~$\mathfrak B_1$. For any~$x \in E,$
\begin{align*}
&\mathds E_{X' \sim \nu}\left[ (p_{i,j}^{(1)})^2(x,X') \right]= \int_z \left(I^{(1)}_{i,j}(x,z)-\mathds E_{\tilde X\sim \pi}[I^{(1)}_{i,j}(x,\tilde X)]\right)^2d\nu(z)\\
&= \int_z \big(\int_w h_{i,j}(x,w)P(z,dw)-\int_w  h_{i,j}(x,w)\underbrace{\int_a P(a,dw)d\pi(a)}_{=d\pi(w)}\big)^2d\nu(z)\\
&=\mathds E_{X' \sim \nu}\left[ \mathds E_{X \sim P(X',\cdot)} h_{i,j}(x,X)-\mathds E_{\pi} [h_{i,j}] \right]^2,
\end{align*}
and for any~$y\in E$,
\begin{align*}
&\mathds E_{\tilde X \sim \pi}\left[ (p_{i,j}^{(1)})^2(\tilde X,y) \right]= \int_x \left(I^{(1)}_{i,j}(x,y)-\mathds E_{\tilde X\sim \pi}[I^{(1)}_{i,j}(x,\tilde X)]\right)^2d\pi(x)\\
&= \int_x \big(\int_w h_{i,j}(x,w)P(y,dw)-\int_w  h_{i,j}(x,w)\underbrace{\int_a P(a,dw)d\pi(a)}_{=d\pi(w)}\big)^2d\pi(x)\\
&=\mathds E_{\tilde X \sim \pi}\left[ \mathds E_{X \sim P(y,\cdot)} h_{i,j}(\tilde X, X)-\mathds E_{  \pi} [h_{i,j}] \right]^2.
\end{align*}
%where~$ \mathrm{Var}_{X' \sim \mathds Q_1}$ denotes the variance with respect to the measure~$\mathds Q_1$ defined by ~$\forall A \in \Sigma,\quad \mathds Q_1(A) =  \int_z P(z,A)d\nu(z)$. 
%From Assumption~\ref{assumption2}.$(ii)$, we have that~$\forall A \in \Sigma,\quad \mathds Q(A) = \nu(A).$ Hence, we get that for any~$x \in E,$
%\begin{align*}&\mathds E_{X' \sim \nu}\left[ (p_{i,j}^{(1)})^2(x,X') \right]\leq \mathrm{Var}_{X' \sim \mathds \nu}(h_{i,j}(x,X')) = \mathds E_{X' \sim \nu}[p_{i,j}^2(x,X')].\end{align*}
%Noticing that Assumption~\ref{assumption2}.$(ii)$ gives that~$\pi=\nu$, 
Hence we get that \begin{align}\label{eq:desc-induction} \mathfrak B^2_1&:=\max\left[ \max_{i} \left\|\sum_{j=i+1}^n \mathds E_{X \sim \nu}\left[ (p_{i,j}^{(1)})^2(\cdot,X) \right] \right\|_{\infty} , \; \max_{j} \left\| \sum_{i=1}^{j-1} \mathds{E}_{X \sim \pi}[(p_{i,j}^{(1)})^2(X,\cdot)] \right\|_{\infty}\right]\\
&\leq  B_n^2 ,\notag \end{align}
where we recall that
\begin{align*}B_n^2= \max\Bigg[ &\sup_{0 \leq k\leq t_n} \; \max_{i} \;\sup_x \sum_{j=i+1}^n \mathds E_{X' \sim \nu}\left[ \mathds E_{X \sim P^k(X',\cdot)} h_{i,j}(x,X)-\mathds E_{  \pi} [h_{i,j}]\right]^2 , \\
&\sup_{0 \leq k\leq t_n}\; \max_{j}\; \sup_{y}\sum_{i=1}^{j-1}\mathds E_{\tilde X \sim \pi}\left[ \mathds E_{X \sim P^k(y,\cdot)} h_{i,j}(\tilde X, X)-\mathds E_{  \pi} [h_{i,j}]\right]^2 \Bigg].\end{align*}

\noindent
This allows us to get a concentration inequality similar to the one of the previous subsection, namely for any~$u>0$, it holds with probability at least~$1-(1+e)e^{-u}$,
\begin{equation*}
 \sum_{i=1}^{n-2} \sum_{j=i+1}^{n-1} h^{(1)}_{i,j}(X_i,X_{j-1},X_{j})  \leq  \kappa \left( A \sqrt{n}  \sqrt{u} + (A  +  B_n\sqrt n) u +2A\sqrt n u^{3/2} + A  u^2\right)
\end{equation*}
\end{itemize}

\section{Proof of Lemma~\ref{lemma:2a-2c}}

\label{proof:lemma-2a-2c}

Using Assumption~\ref{assumption4}, we have that~$\mathds E_{\pi}[h_{i,j}]= \int_{x_i} P^{i-j+t_n}(X_{j-t_n},dx_i)  \int_{x_j} h_{i,j}(x_i,x_j) d\pi(x_j)$. Hence we get,
\begin{align*}
(2a)&:= \left| \sum_{j=2}^{n} \; \sum_{i=(j-t_n+1)\vee 1}^{j-\floor{\frac{t_n}{2}}}  \mathds E_{j-t_n}\left[ h_{i,j}(X_i,X_j)\right]-\mathds E_{\pi}\left[h_{i,j}\right] \right|\\
&\leq\sum_{j=2}^{n} \; \sum_{i=(j-t_n+1)\vee 1}^{j-\floor{\frac{t_n}{2}}}  \left| \int_{x_i} \int_{x_j}  h_{i,j}(x_i,x_j) P^{i-j+t_n}(X_{j-t_n},dx_i)\left( P^{j-i}(x_i,dx_j) -d\pi(x_j)\right)\right|\\
%&\leq\sum_{j=2}^{n} \; \sum_{i=(j-t_n+1)\vee 1}^{j-\floor{\frac{t_n}{2}}}  \|h_{i,j}\|_{\infty}\int_{x_i}    P^{i-j+t_n}(X_{j-t_n},dx_i) \int_{x_j}\left| P^{j-i}(x_i,dx_j) -d\pi(x_j)\right|\\
&\leq\sum_{j=2}^{n} \; \sum_{i=(j-t_n+1)\vee 1}^{j-\floor{\frac{t_n}{2}}}  \|h_{i,j}\|_{\infty} \underbrace{\int_{x_i}    P^{i-j+t_n}(X_{j-t_n},dx_i)}_{=1} \underbrace{\sup_{y}\int_{x_j}\left| P^{j-i}(y,dx_j) -d\pi(x_j)\right|}_{=\sup_y \| P^{j-i}(y,\cdot)-\pi\|_{\mathrm{TV}}}\\
%&\leq\sum_{j=2}^{n} \; \sum_{i=(j-t_n+1)\vee 1}^{j-\floor{\frac{t_n}{2}}}  \|h_{i,j}\|_{\infty}\sup_y \| P^{j-i}(y,\cdot)-\pi\|_{\mathrm{TV}}\\
&\leq\sum_{j=2}^{n} \; \sum_{i=(j-t_n+1)\vee 1}^{j-\floor{\frac{t_n}{2}}} \|h_{i,j}\|_{\infty}L\rho^{j-i}\leq\sum_{j=2}^{n} \;\sum_{i=(j-t_n+1)\vee 1}^{j-\floor{\frac{t_n}{2}}} \|h_{i,j}\|_{\infty}L\rho^{t_n/2}\leq LAt_n,
\end{align*}
where we used that~$\rho^{t_n/2}\leq\rho^{r\log(n)/2} = n^{r \log(\rho)/2}\leq n^{-1}~$. Indeed ~$1+r\log(\rho)/2<0$ because we choose~$r~$ such that~$r>2(\log(1/\rho))^{-1}.$ With an analogous approach, we bound the term $(2c)$ as follows.
\begin{align*}
(2c)&:= \left| \sum_{j=2}^{n} \; \sum_{i=(j-\floor{\frac{t_n}{2}}+1)\vee 1}^{j-1}  \mathds E_{j-t_n}\left[ h_{i,j}(X_i,X_j)\right]-\mathds E\left[h_{i,j}( X_i, X_j)]\right] \right|\\
&\leq\sum_{j=2}^{n} \; \sum_{i=(j-\floor{\frac{t_n}{2}}+1)\vee 1}^{j-1}   \left|\int_{x_j} \int_{x_i} P^{j-i}(x_i,dx_j) h_{i,j}(x_i,x_j) \left(P^{i-j+t_n}(X_{j-t_n},dx_i) -d\pi(x_i)\right)\right|\\
&\leq\sum_{j=2}^{n} \; \sum_{i=(j-\floor{\frac{t_n}{2}}+1)\vee 1}^{j-1} \|h_{i,j}\|_{\infty}  \sup_{z} \|P^{i-j+t_n}(z,\cdot) -\pi\|_{\mathrm{TV}}\\
&\leq\sum_{j=\floor{\frac{t_n}{2}}}^{n} \; \sum_{i=(j-\floor{\frac{t_n}{2}}+1)}^{j-1} \|h_{i,j}\|_{\infty}  L\rho^{i-j+t_n} +\sum_{j=2}^{\floor{\frac{t_n}{2}}} \; \sum_{i=(j-\floor{\frac{t_n}{2}}+1)\vee 1}^{j-1} \|h_{i,j}\|_{\infty}  L\rho^{i-j+t_n}  \\
&\leq\sum_{j=\floor{\frac{t_n}{2}}}^{n} \; \sum_{i=(j-\floor{\frac{t_n}{2}}+1)}^{j-1} \|h_{i,j}\|_{\infty}  L\rho^{t_n/2} +t_n^2\|h_{i,j}\|_{\infty}  L  \leq LA\left( t_n^2 + nt_n \rho^{t_n/2}\right)\leq 2LA t_n^2,
\end{align*}
where we used that~$\rho^{t_n/2}\leq\rho^{r\log(n)/2} = n^{r \log(\rho)/2}\leq n^{-1}~$.

\section{Proof of Proposition~\ref{prop:tau}} \label{proof:prop-regeneration} Since the split chain has the same distribution as the original Markov chain, we get that~$(\widetilde X_i)_i$ is~$\psi$-irreducible for some measure~$\psi$ and uniformly ergodic. From~\cite[Theorem~16.0.2]{tweedie}, Assumption~\ref{assumption1} ensures that for every measurable set~$A\subset E\times \{0,1\}$ such that~$\psi(A)>0$, there exists some~$\kappa_A>1$ such that \[\sup_{x} \mathds E[\kappa_A^{\tau_A}|\widetilde X_1=x]<\infty,\]where~$\tau_A:=\inf \{n\geq 1 \; : \; \widetilde X_n \in A\}$ is the first hitting time of the set~$A$. Let us recall that~$T_1$ and~$T_2$ are defined as hitting times of the atom of the split chain~$E\times \{1\}$ which is accessible (i.e. the atom has a positive~$\psi$-measure). Hence, there exist~$C>0$ and~$\kappa>1$ such that, \[\sup_{x}\mathds E[\kappa^{\tau_{E\times\{1\}}}|\widetilde X_1=x]=\sup_{x} \mathds E[\exp(\tau_{E\times\{1\}}\log(\kappa))|\widetilde X_1=x]\leq C.\]
Considering~$k\geq 1$ such that~$C^{1/k} \leq 2$, a straight forward application of Jensen inequality gives that~$\max(\|T_1\|_{\psi_1},\|T_2\|_{\psi_1})\leq k/\log(\kappa)$.

\end{document}